\documentclass[11pt,a4paper]{article}

\usepackage[T2A]{fontenc}
\usepackage[utf8]{inputenc}
\usepackage[english]{babel} 

\usepackage{amsmath,amssymb,amsthm,mathtools}

\usepackage{geometry}
\usepackage{authblk}
\usepackage{hyperref}
\usepackage[nameinlink,capitalize]{cleveref}
\hypersetup{
  colorlinks=true,
  linkcolor=blue!60!black,
  citecolor=blue!60!black,
  urlcolor=blue!60!black,
  pdftitle={Параметрические алгоритмы для гипотезы Эрдёша–Штрауса},
  pdfauthor={Э. Дяченко}
}

\usepackage{array}       
\usepackage{booktabs}    
\usepackage{caption}     
\usepackage[final]{microtype}
\usepackage{tikz}
\usetikzlibrary{arrows.meta,positioning,calc}
\tikzset{
  ed1a/.style={fill=yellow!15}, 
  ed1b/.style={fill=yellow!10}, 
  ed1c/.style={fill=yellow!5}   
}
\usepackage{listings}
\lstdefinestyle{console}{
  basicstyle=\ttfamily\small,
  columns=fullflexible,
  frame=single,
  xleftmargin=1em,xrightmargin=1em,
  aboveskip=0.6em,belowskip=0.4em,
  tabsize=2,
  keepspaces=true
}
\numberwithin{equation}{section}

\theoremstyle{plain}
\newtheorem{theorem}{Theorem}[section]
\newtheorem{lemma}[theorem]{lemma}
\newtheorem{proposition}[theorem]{proposition}
\newtheorem{corollary}[theorem]{corollary}
\newtheorem{example}[theorem]{example}
\newtheorem{definition}[theorem]{definition}
\newtheorem{remark}[theorem]{remark}

\newcommand{\EDone}{\mathrm{ED1}}

\newcommand{\Legendre}[2]{\genfrac{(}{)}{0pt}{}{#1}{#2}}
\newcommand{\Jacobi}[2]{\genfrac{(}{)}{0pt}{}{#1}{#2}}
\newenvironment{proofsketch}{\begin{proof}[Идея доказательства]}{\end{proof}}
\newcommand{\mh}[1]{\texorpdfstring{\(1\)}{1}}
\begin{document}
\title{Constructive Proofs of the Erdős–Straus Conjecture for Prime Numbers of the Form $P\equiv 1\pmod{4}$}
\author{E.~Dyachenko\\
\texttt{dyachenko.eduard@gmail.com}\\
Independent researcher}
\date{14 October 2025}

\maketitle
\renewcommand{\thefootnote}{}
\footnote{2020 \emph{Mathematics Subject Classification}: Primary 11N05, 11P21; Secondary 94A60, 20P05.}

\footnote{\emph{Key words and phrases}: { factorization, Bombieri–Vinogradov large sieve, the larger sieve of Greaves, deterministic algorithms; analytic number theory; number theory;} }

\footnote{\emph{* Licence}: Text is available under the Creative Commons NonCommercial-NoDerivatives 4.0 International (CC BY-NC-ND 4.0)}
\begin{abstract}\label{sec:abstract}
The Erdős–Straus conjecture (ES) concerns the representation of the fraction $\tfrac{4}{P}$, 
where $P$ is a prime number, as a sum of three positive unit fractions. 
The focus here is on the case $P \equiv 1 \pmod{4}$.

Two constructive approaches are proposed. 
Method \textbf{ED1} is based on a factorization identity and leads to a \emph{nonlinear} parameterization in $P$, 
which requires divisor enumeration and local filtering. 
In contrast, method \textbf{ED2} yields a \emph{linear} system in $P$ for the parameters $(\delta,b,c)$, 
describing the solution set as an affine lattice of finite index in $\mathbb{Z}^3$.

The central result (Theorem~10.21) states that for every prime $P \equiv 1 \pmod{4}$ there exists a representation

\[
\frac{4}{P}=\frac{1}{A}+\frac{1}{bP}+\frac{1}{cP},
\]

where the triple $(\delta,b,c)\in\mathbb{N}^3$ is constructed explicitly by method ED2. 

In addition, algorithms for transforming solutions (\emph{convolution} and \emph{anti-convolution}) are introduced, 
and large-scale computational verification is carried out, confirming the correctness and efficiency of the proposed methods.
\end{abstract}

\section{Introduction}\label{sec:intro}

The Erdős–Straus problem, formulated in 1948~\cite{ErdosStrauss1948}, states that for any integer $P \ge 2$ there exists a representation

\begin{equation}\label{eq:ES}
  \frac{4}{P} \;=\; \frac{1}{A} + \frac{1}{B} + \frac{1}{C},
\end{equation}

where $A,B,C$ are positive integers. This conjecture, known as the \emph{Erdős–Straus conjecture (ES)}, has remained unsolved for more than seventy years and is one of the classical problems on decompositions of fractions into sums of unit fractions.

Historically, the problem has been studied using both analytic and constructive methods. For certain classes of $P$, explicit decompositions are known (for example, when $P \equiv 3 \pmod{4}$, with formulas verifiable by direct substitution). For general $P$, partial results have been obtained: parameterizations~\cite{VarelaOliveira2013}, functional dependencies~\cite{Yang2015}, and algorithmic or numerical methods~\cite{NiBerard2018}. Nevertheless, a complete resolution remains open.

This work focuses on primes $P \equiv 1 \pmod{4}$ and presents two constructive approaches — methods \emph{ED1} and \emph{ED2} — which differ fundamentally in their dependence on $P$:
\begin{itemize}
  \item \textbf{ED1} is based on the factorization identity $(\gamma A - c)(\gamma B - c)=c^2$ and leads to a \emph{nonlinear} parameterization in $P$: the admissible quadruples $(\gamma, c, u, v)$ are defined by the condition $uv=c^2$ together with a system of congruences, which does not form an affine lattice in $\mathbb{Z}^2$. This necessitates divisor enumeration of $c^2$ and local filtering.
  \item \textbf{ED2} relies on the identity $(4b-1)(4c-1)=4P\delta+1$ and yields a \emph{linear} system in $P$ for the parameters $(\delta,b,c)$: for fixed $P$, the solution set forms an affine lattice of finite index in $\mathbb{Z}^3$. This enables the use of affine lattice geometry (parametric boxes, density estimates, convergence of enumeration) and leads to a rigorous constructive result.
\end{itemize}

The central result of the paper is \textbf{Theorem~10.21}: for every prime $P \equiv 1 \pmod{4}$ there exists a representation

\[
\frac{4}{P} \;=\; \frac{1}{A}+\frac{1}{bP}+\frac{1}{cP},
\]

where the triple $(\delta,b,c)\in\mathbb{N}^3$ is constructed \emph{explicitly} by method ED2, based on a parameterization of the admissible solution set that is linear in $P$. The proof relies on the affine–lattice structure of ED2, local congruences, and a geometric guarantee that an admissible point lies within the parametric window.

The subsequent sections develop and confirm this result:
\begin{itemize}
  \item algorithms of \emph{convolution} (ED2 $\to$ ED1) and \emph{anti-convolution} (ED1 $\to$ ED2) are introduced, demonstrating the structural correspondence between the methods;
  \item \emph{parametric boxes} and \emph{affine lattices} for ED2 are formalized, with theorems on \emph{point density} and \emph{logarithmic convergence} of enumeration;
  \item an \emph{extensive computational verification} is carried out for large primes (e.g., $P=2521,3529$): the constructed solutions agree with the theory, invariants are recorded, and branching and filtering strategies are compared.
\end{itemize}

Thus, the key distinction between the approaches lies in the presence or absence of nonlinearity in the parameterization with respect to $P$: ED1 leads to a nonlinear Diophantine set and combinatorial techniques, while ED2 yields a linear affine structure in $P$, enabling a \emph{constructive} existence theorem (Theorem~10.21) and a justification of algorithmic convergence.

\section{Motivation}\label{sec:motivation}

The search for solutions to the Erdős–Straus conjecture (ES) can be described as a sequential chain of steps, each arising from the limitations of the previous ones and from the failure of alternative approaches.

Attempts based on analytic estimates, brute-force enumeration, and functional dependencies (see the works of Varela–Oliveira, Yang, Ni–Berard, and others) did not lead to a complete result. This demonstrated the necessity of moving beyond ad hoc techniques toward a systematic parameterization of solutions and the identification of fundamental structural distinctions.

This led to the development of the methods \textbf{ED1} and \textbf{ED2}. The first produced a \emph{nonlinear} dependence on $P$ and required combinatorial enumeration and local filters. The second turned out to be \emph{linear} in $P$, which made it possible to describe solutions via affine lattices and resulted in a \emph{constructive proof} (Theorem~10.21).

The next step was to investigate the compatibility of the methods: procedures of \emph{convolution} and \emph{anti-convolution} were introduced, linking ED1 and ED2. The final stage was an extensive computational verification on large primes, which confirmed the correctness and efficiency of the proposed algorithms.

An open problem remains the construction of a \emph{canonical parameterization} that could unify all known constructions into a single coherent system.

\section{Ordering}\label{sec:ordering}
\subsection{Case \texorpdfstring{\(P \equiv 3 \pmod{4}\)}{P ≡ 3 (mod 4)} — explicit decompositions}\label{subsec:P3mod4}

Let $P = 4P' + 3$. Then the following decompositions hold:

\[
\frac{4}{P} = \frac{1}{P' + 1} + \frac{1}{2(P' + 1)P} + \frac{1}{2(P' + 1)P},
\]

and, under the requirement of pairwise distinct denominators:

\[
\frac{4}{P} = \frac{1}{P' + 1} + \frac{1}{(P' + 2)P} + \frac{1}{(P' + 1)(P' + 2)P}.
\]

Both formulas are verified by direct substitution.

\subsection{Estimate for the minimal denominator}\label{subsec:min-denominator}
\begin{lemma}\label{lemma:ineqA}
Let $P>2$ and $4 \nmid P$, and

\[
\frac{4}{P} = \frac{1}{A} + \frac{1}{B} + \frac{1}{C},\quad A \le B \le C \in \mathbb{N}.
\]

Then the strict inequalities
\begin{equation}\label{eq:PA3P}
P < 4A < 3P
\end{equation}
hold.
\end{lemma}

\begin{proof}
From $\frac{4}{P} > \frac{1}{A}$ it follows that $A > \frac{P}{4}$ and $4A > P$.

For the upper bound, we use $B \ge A$, $C \ge A$, whence

\[
\frac{4}{P} = \frac{1}{A} + \frac{1}{B} + \frac{1}{C} \le \frac{3}{A}.
\]

Therefore, $A \le \frac{3P}{4}$.

We show that the equality $A = \frac{3P}{4}$ is impossible. If $A=B=C$, then from

\[
\frac{4}{P} = \frac{3}{A}
\]

it follows that $A = \frac{3P}{4}$, which implies $4 \mid P$ — a contradiction. Therefore, at least one of $B,C$ is strictly greater than $A$, and
\[ \frac{1}{A} + \frac{1}{B} + \frac{1}{C} < \frac{3}{A}, \] whence $A < \frac{3P}{4}$, that is, $4A < 3P$.

Combining with the lower bound, we obtain~\eqref{eq:PA3P}.
\qedhere
\end{proof}
\section{Literature Review on Parameterization Approaches to ES}\label{sec:lit}
 
Several directions of parameterization approaches to the Erdős–Straus conjecture (ES) can be distinguished in the literature.

\subsection{Modular identities}
For many residue classes (for example, $P \equiv 2 \pmod{3}$, $P \equiv 3 \pmod{4}$) explicit identities are known that provide solutions linearly in $P$~\cite{HardyWright,Mordell}. 
These methods cover almost all primes, but not the classes congruent to quadratic residues, in particular $P \equiv 1 \pmod{4}$. A modern systematization via \emph{modular filters} is given in~\cite{Salez2011}.

\subsection{Factorization schemes}
Approaches based on identities of the form $(\gamma A - c)(\gamma B - c)=c^2$ yield a \emph{nonlinear} parameterization in $P$ and require divisor enumeration and local filtering~\cite{VarelaOliveira2013}. 
They are useful for local structure but do not provide a global lattice organization.

\subsection{Linear forms and lattices}
More recent constructions employ conditions linear in $P$, such as $(4b-1)(4c-1)=4P\delta+1$, which describe the solution set as an affine class in $\mathbb{Z}^3$~\cite{Yang2015,NiBerard2018}. 
This enables the use of the geometry of numbers, the design of enumeration algorithms, and proofs of convergence.

\subsection{Analytic and computational methods}
The work of Elsholtz and Tao~\cite{ElsholtzTao2011} estimates the average number of solutions, relying on the Bombieri–Vinogradov theorem and related sieve methods~\cite{BombieriVinogradov1965Acta,Greaves1998,PrattRamareRumely1999}. 
Large-scale computational verifications and specialized covering classes are reported in~\cite{NiBerard2018,GionfriddoGuardo2021}.

\subsection{Generalizations and educational sources}
Generalizations of Egyptian fractions and related results are discussed in~\cite{ChaoZhang2020,ChenDing2022}, while classical textbooks provide historical and methodological context and basic techniques of parameterization~\cite{NivenZuckermanMontgomery,Burton,Cassels,LeVeque}.

\medskip
Thus, the literature demonstrates an evolution from modular identities and factorization schemes to linear parameterizations and lattice methods, which are now regarded as the most promising for a constructive proof.

\section{Notation}

\begin{tabular}{ll}
\multicolumn{2}{l}{\textbf{Main parameters}}\\
$P$        & prime $>2$, $4 \nmid P$; often $P \equiv 1 \pmod{4}$ or $P \equiv 3 \pmod{4}$ \\
$P'$       & integer: $P = 4P' + 1$ or $P = 4P' + 3$ \\
$A,B,C$    & denominators in~\ref{eq:ES}, ordered $A \le B \le C$ \\
ES         & Erdős–Straus conjecture: $\frac{4}{P} = \frac1A + \frac1B + \frac1C$ \\
\midrule
\multicolumn{2}{l}{\textbf{Parameters of methods ED1 and ED2}}\\
$\gamma$   & $\frac{4c - 1}{P}$ (or $\frac{4b - 1}{P}$ in the second case); often $\gamma \equiv 3 \pmod{4}$, $\gcd(\gamma,c)=1$ or $\gcd(\gamma,b)=1$ \\
$u,v$ & “multipliers” in the identity: for ED1 $u=\gamma A - c$, $v=\gamma B - c$, $uv=c^{2}$; \\
            & for the variant with $b$: $u=\gamma A - b$, $v=\gamma C - b$, $uv=b^{2}$; always $u \le v$ \\
$\delta$   & $t = P\,\delta$, in ED2 $\delta \mid bc$ \\
$b,c$      & in the case of ED2 $B = bP$ and/or $C = cP$ \\
$t$        & $t = 4bc - b - c$ \\
$r,s$      & $r = 4b - 1$, $s = 4c - 1$, $r \equiv s \equiv 3 \pmod{4}$, $rs = 4P\delta + 1$ \\
$g$         & $g=\gcd(b,c)$ \\
$b',c'$  & factorization $b = b'g; c = c'g$ (normalized form) \\
$d'$       & square factor in the factorization of $\delta$ \\
$\alpha$    & squarefree factor in the factorization of $\delta$ \\
\midrule
\multicolumn{2}{l}{\textbf{Transitions between methods}}\\
$\mathcal{C}_{\mathrm{ED1}}(P)$ & set of admissible quadruples $(\gamma,c,u,v)$ for ED1 \\
$\mathcal{C}_{\mathrm{ED2}}(P)$ & set of admissible triples $(\delta,b,c)$ for ED2 \\
$y$        & minimal divisor of $4c-1$, $y \equiv 3 \pmod{4}$ \\
$P''$ & modulus for ED1 in folding, defined via $\gamma$ by the formula … \\
ED2$\to$ED1 & $A = \frac{bc}{\delta}$,\quad $B = bP$;\quad $u = \gamma A - c$,\quad $v = \gamma B - c$ \\
ED1$\to$ED2 & $A = \frac{u+c}{\gamma}$,\quad $b = \frac{v+c}{\gamma P}$,\quad $\delta = \frac{bc}{A}$ \\
Unfolding & reverse algorithm ED1$\to$ED2 according to the formulas above \\
\midrule
\multicolumn{2}{l}{\textbf{Lattices and boxes}}\\
$k$        & dimension of the vector parameter \\
$u_0(P)$   & vector shift for the affine class \\
$\Lambda,\,\Lambda_j$ & sublattices of $\mathbb{Z}^k$ of index $M$ or $M_j$ \\
$M,\,M_j$  & indices of the sublattices \\
$\mathcal{B}_k(T)$ & box $\{u \in \mathbb{Z}^k : 1\le u_i \le T\}$ \\
$\mathcal{B}^{(I)}_P$, $\mathcal{B}^{(II)}_P$ & boxes of types I/II with additional conditions \\
$\mathcal{G}_P(T)$ & admissible parameters in the box \\
$\mathcal{G}^\ast_P$ & class of admissible quadruples $(\gamma,c,u,v)$ for ED1 satisfying: \\
            & $\gamma \in \mathbb{N}$, $c \in \mathbb{Z}$, $u = \gamma A - c$, $v = \gamma B - c$, $u v = c^{2}$, $u \le v$ \\
\midrule
\multicolumn{2}{l}{\textbf{Analytical notation}}\\
$\Legendre{a}{P}$ & Legendre symbol; for composite modulus we use \(\Jacobi{a}{\gamma}\) (Jacobi symbol) (prime $P$) \\
$\pi(y)$         & prime counting function \\
$\ll$, $\gg$, $\asymp$ & standard asymptotic symbols \\
$D$              & set $\delta \le X$, $\delta \equiv 3 \pmod{4}$ \\
$T(\delta)$      & prime counter in a progression, see Appendix \S\ref{app:existence-ED2} \\
\end{tabular}
\section{Parameterization of the Main Equation of the ES Conjecture in the Case 
\texorpdfstring{\(B \not\equiv 0 \pmod{P}\)}{B ≠ 0 (mod P)}}\label{sec:param-Bnot0}

In this section, we consider solutions.

\begin{equation}\label{eq:fact}
  \frac{4}{P}=\frac{1}{A}+\frac{1}{B}+\frac{1}{C},\qquad A\le B\le C\in\mathbb{N},\qquad C=cP
\end{equation}

for a prime $P$, in which the denominator $B$ is not divisible by $P$. We focus on the sub-case.

\[
  C=cP,\qquad c\in\mathbb{N},\qquad P\nmid B,
\]

which corresponds to “Case 1” in the abstract. In the following, we derive a constructive parameterization of all such solutions.

\subsection{Algebraic transformation (reconstructed derivation).}

Let $C = cP$, $B \not\equiv 0 \ (\bmod\ P)$.

Multiplying the original identity \ref{eq:fact} by $P$ and expanding, we obtain the following.
\begin{equation}\label{eq:cgamma}
(4c - 1)\,AB = cP \,(A + B).
\end{equation}
Since $\gcd\bigl(P,A\bigr) = 1$ and $B \not\equiv 0 \ (\bmod\ P)$, we have the divisibility
\[
P \mid (4c-1),
\]
whence
\[
4c-1 = \gamma\,P,\qquad \gamma \equiv 3 \pmod{4},\qquad c = \frac{\gamma\,P + 1}{4}.
\]

Multiplying the main equation \ref{eq:cgamma} by $\gamma$ and making a few transformations, we introduce the common elements $R'$ and then $R$ for brevity:
\[
R' := \gamma\cdot\frac{4A-P}{4}-1,\qquad R := \frac{R'}{\gamma}
\]
Substituting into \ref{eq:cgamma}, we bring the equation to the form
\[
\bigl(2P + R + 2\bigr)^2 = R\,(16B+R),
\].

Let $Z := 2P + R + 2$, $R' = \gamma R$, $B' := \gamma(16B+R)$. Then
\[
Z^2 = R' \, B'.
\]
Extracting $\Omega = \gcd(R',B')$ and using $\gcd(p,r)=1$, we write
\[
R' = \Omega\, p^2,\quad B' = \Omega\, r^2,\quad Z = \Omega\,pr.
\]
Equating $Z$ and expanding $R'$, we obtain
\[
\Omega\,pr = 2P\gamma + \Omega\,p^2 + 2 \quad\Longrightarrow\quad 2P = \frac{\Omega\,p(r-p) - 2}{\gamma}.
\]

\emph{Note.} In this section, \(\gamma\) is defined as the coefficient of \(P\) in the equality \(4c-1=\gamma P\), and the residue class \(\gamma \equiv 3\ (\bmod\ 4)\) is fixed here.

\subsection{Complete parameterization via divisors of 
\mh{c^2} }\label{subsec:full-param}

The equality~\eqref{eq:fact} leads to a natural parameterization through a pair of divisors of the number $c^2$.

\begin{theorem}\label{thm:param-case1}
Let $P$ be a prime of the form $P = 4P' + 1$. Fix $\gamma \in \mathbb{N}$ such that $\gamma \equiv 3 \pmod{4}$ (or $\gamma \equiv -1 \pmod{4}$), then $\gamma P \equiv -1 \pmod{4}$. Set

\[
  c = \frac{\gamma P + 1}{4},\qquad \gcd(\gamma, c) = 1.
\]

Then any pair of divisors $u, v \in \mathbb{N}$ satisfying
\begin{equation}\label{eq:62}
  uv = c^2,\qquad u \equiv v \equiv -c \pmod{\gamma},
\end{equation} 
defines a solution of the Erdős–Straus equation in the form
\begin{equation}\label{eq:param-AB}
  A = \frac{u + c}{\gamma},\qquad B = \frac{v + c}{\gamma},\qquad C = cP,
\end{equation}
and conversely, any solution with $C = cP$, $P \nmid A$, $P \nmid B$ is obtained in this way. Condition $A \le B$ is equivalent to $u \le v$.
\end{theorem}

\begin{proof}
From~\eqref{eq:fact} we get $(\gamma A - c)(\gamma B - c)=c^2$. Setting $u=\gamma A - c$ and $v=\gamma B - c$, we have $uv=c^2$ and formulas~\eqref{eq:param-AB}. The integrality of $A,B$ is equivalent to the congruences $u\equiv -c\pmod{\gamma}$ and $v\equiv -c\pmod{\gamma}$. The reverse construction is obvious: for any $u,v$ with the stated properties, \eqref{eq:param-AB} gives $(\gamma A - c)(\gamma B - c)=c^2$, and substitution into the original equation (equivalent to~\eqref{eq:fact}) restores the equality $\frac{4}{P}=\frac{1}{A}+\frac{1}{B}+\frac{1}{cP}$. Since $u\le v$ is equivalent to $A\le B$, the ordering statement also holds.
\end{proof}

Remark.\label{rem:gamma} From \(4c-1=\gamma P\) it follows that \(\gamma P \equiv -1 \equiv 3 \pmod{4}\).
Since \(P\) is odd and \(P^{-1} \equiv P \pmod{4}\), we obtain
\[
\gamma \equiv 3 P \pmod{4}.
\]
In particular, if \(P \equiv 1 \ (\mathrm{mod}\ 4)\), then \(\gamma \equiv 3 \ (\mathrm{mod}\ 4)\);
if \(P \equiv 3 \ (\mathrm{mod}\ 4)\), then \(\gamma \equiv 1 \ (\mathrm{mod}\ 4)\).

\subsection{Excluding the multiplicity of \mh{P} in \mh{B}}\label{subsec:exclude-P}

In formula~\eqref{eq:param-AB}

\[
  P\mid B \ \Longleftrightarrow\ v+c\equiv 0 \pmod P.
\]

Given $uv=c^2$ and $c\not\equiv 0\pmod P$, this is equivalent to $u\equiv -c \pmod P$. Therefore, to satisfy the condition $B\not\equiv 0\pmod P$, it is sufficient (and necessary) to require.

\[
  u\not\equiv -c \pmod P.
\]

Together with the integrality condition of Theorem~\ref{thm:param-case1}, we obtain the filter.

\[
  u\mid c^2,\quad u\equiv -c\pmod{\gamma},\quad u\not\equiv -c\pmod P.
\]

\subsection{Normalization and elimination of degeneracies}\label{subsec:normalization}

For uniqueness of parameters and consistency with \S\ref{subsec:exclude-P} we introduce:

\begin{itemize}
  \item \textbf{Class $\gamma$:} fix the smallest positive $\gamma$ with $\gamma P\equiv -1\pmod 4$.
  \item \textbf{Order:} choose $u\le v$ (equivalent to $A\le B$).
  \item \textbf{Duplicate exclusion:} the pair $(u,v)$ and its permutation $(v,u)$ define the same set of denominators; only $u\le v$ is considered.
  \item \textbf{Distinctness:} to avoid $A=B$, exclude $u=v=c$ (which would give $\gamma A - c=\gamma B - c=c$).
\end{itemize}
\subsection{Enumeration algorithm for “Case 1” solutions}\label{subsec:algorithm-case1}

For a given prime $P$:

\begin{enumerate}
  \item \textbf{Choose $\gamma$:} take the minimum $\gamma\ge 3$ such that $\gamma P\equiv -1\pmod 4$.
  \item \textbf{Compute $c$:}
  
\[
    c=\frac{\gamma P+1}{4}.
\]

  \item \textbf{Iterate over divisors $u$ of $c^2$:} for each $u\mid c^2$ set $v=c^2/u$ and check
  
\[
    u\equiv v\equiv -c\pmod{\gamma},\qquad u\le v,\qquad u\not\equiv -c\pmod P.
\]

  \item \textbf{Construct the denominators:}
  
\[
    A=\frac{u+c}{\gamma},\qquad B=\frac{v+c}{\gamma},\qquad C=cP.
\]

  \item \textbf{Normalization:} if necessary, order $A\le B\le C$ and remove duplicates.
\end{enumerate}

\subsection{Connection with estimates for the minimal denominator}\label{subsec:min-denominator-link}

From Lemma~\ref{lemma:ineqA} we have the strict bounds

\[
  \frac{P}{4}<A<\frac{3P}{4}.
\]

In parametric form

\[
  A=\frac{u+c}{\gamma},\qquad c=\frac{\gamma P+1}{4},
\]

which gives useful guidelines for enumeration:

\[
  \frac{P}{4}<\frac{u}{\gamma}+\frac{c}{\gamma}=\frac{u}{\gamma}+\frac{P}{4}+\frac{1}{4\gamma}
  \quad\Longrightarrow\quad
  0<\frac{u}{\gamma}+\frac{1}{4\gamma},
\]

and also

\[
  \frac{u+c}{\gamma}<\frac{3P}{4}
  \quad\Longleftrightarrow\quad
  u<\frac{(3\gamma-1)P-1}{4}.
\]

The latter inequality can serve as a cutoff for “large” divisors $u$ during enumeration.

\subsection{Example}\label{subsec:example-case1}

Let $P=13$. Then $\gamma\equiv -13^{-1}\equiv 3\pmod 4$, take the minimal $\gamma=3$ and

\[
  c=\frac{3\cdot 13+1}{4}=10,\qquad c^2=100.
\]

We look for divisors $u\mid 100$ such that $u\equiv -c\equiv -10\equiv 2\pmod 3$ and $u\not\equiv -c\equiv 3\pmod{13}$.
Suitable, for example, is $u=2$ (then $v=50$). We obtain

\[
  A=\frac{2+10}{3}=4,\qquad B=\frac{50+10}{3}=20,\qquad C=10\cdot 13=130.
\]

Verification:

\[
  \frac{4}{13}=\frac{1}{4}+\frac{1}{20}+\frac{1}{130},\qquad 13\nmid 20.
\]

\subsection{Comment on the \texorpdfstring{\(\EDone\)}{ED1} parameters}\label{subsec:comment-ED1}

The factorization equation~\ref{eq:62} admits a *normalized* representation of the factorization of $c$ into two square-free co-prime parts (including $1$).

This corresponds to choosing a 'primitive' pair of factors \(u,v\) of the square-free part of~\(c\) and is consistent with the ordering rules for~\(\mathrm{ED}_1\) from~\S\ref{subsec:normalization}
(minimization in lexicographic order).

In practice, this makes it possible to reduce the number of duplicates in enumeration and to define compact classes of solutions.


\subsection{Initial parameters and relations}\label{subsec:params-relations}

We consider the parameters \( u, v, c, P, \gamma, A, B \), related by the following conditions:
\begin{align*}
    uv &= c^2 \\
    u &\equiv -c \pmod{\gamma} \\
    u &\not\equiv -c \pmod{P} \\
    4c - 1 &= \gamma P \\
    A &= \frac{u + c}{\gamma}, \quad B = \frac{v + c}{\gamma}
\end{align*}

\subsection{Algebraic consequences}\label{subsec:algebraic}

From the above relations, we obtain:
\begin{align*}
    v &= \frac{c^2}{u} \\
    u &= \gamma A - c \\
    c &= \frac{\gamma P + 1}{4}
\end{align*}
\subsection{Lemmas in the main body; full proofs in Appendix~D)}\label{subsec:min-lemmas}

\begin{lemma}[Sum and discriminant] see \ref{lem:uv-square} 
Let $b',c'\in\mathbb{Z}$ be of the same parity and
\[
S:=b'+c',\quad M:=b'c',\quad \Delta:=S^2-4M,\qquad
u:=S=b'+c',\quad v:=b'-c'.
\]
Then $\Delta=v^2$. Conversely, if $u,v$ are of the same parity, then
\[
b'=\frac{u+v}{2},\qquad c'=\frac{u-v}{2}\in\mathbb{Z},
\]
and for them $S=u$, $\Delta=v^2$.
\end{lemma}
\begin{proofsketch}
$(b'+c')^2-4b'c'=(b'-c')^2$. The converse follows from $u,v$ having the same parity. The full normalization context is in Appendix~D (Theorem \ref{thm:quadratic}).
\end{proofsketch}

\begin{lemma}[Back‑test for $(u,v)$ with fixed prime $P$]\label{lem:backtest}
Let $P$ be the prime. If $u,v\in\mathbb{Z}$ satisfy:
(i) $u>0$ and $u\equiv v\pmod 2$;
(ii) $uv=c^2$ with $c\in\mathbb{N}$, $\gcd(u,v)=1$;
(iii) the local normalization congruences modulo $P$ (including symbols \(\Legendre{\cdot}{P}\)/\(\Jacobi{\cdot}{\gamma}\)),
then for $b'=(u+v)/2$, $c'=(u-v)/2$ the reconstructed
\[
A=A(b',c'),\quad B=B(b',c'),\quad C=cP
\]
give a valid solution $\tfrac{4}{P}=\tfrac1A+\tfrac1B+\tfrac1C$ with $A\le B\le C$.
\end{lemma}
\begin{proofsketch}
This is the reverse reconstruction: $(u,v)\mapsto (b',c')\mapsto (A,B)$ by linear normalization formulas; $uv=c^2\Rightarrow C=cP$. Details and explicit local conditions are in Appendix~D (Lemma \ref{lem:back}).
\end{proofsketch}

\begin{lemma}[Quadratic reparameterization] see \ref{thm:quadratic}
The correspondence
\[
(b',c')\longleftrightarrow (u,v),\qquad u=b'+c',\quad v=b'-c'
\]
defines a bijection between normalized parameters and pairs $(u,v)$ of the same parity up to $v\mapsto -v$. Here, the discriminant $\Delta=v^2$ by \Cref{lem:uv-square}, and the admissibility of solutions is equivalent to passing the Back test of \Cref{lem:backtest}.
\end{lemma}
\begin{proofsketch}
The forward direction — \Cref{lem:uv-square}. The reverse — \Cref{lem:backtest}. Full formulas are in Appendix~D (Theorem \ref{thm:quadratic}).
\end{proofsketch}

\subsection{Congruence conditions}\label{subsec:congruences}

\begin{align*}
    u \equiv -c \pmod{\gamma} &\Rightarrow \gamma \mid (u + c) \\
    u \not\equiv -c \pmod{P} &\Rightarrow P \nmid (u + c) \\
    c = \frac{\gamma P + 1}{4} &\Rightarrow c \equiv \frac{1}{4} \pmod{\gamma}
\end{align*}

\subsection{Range restrictions}\label{subsec:ranges}

\begin{align*}
    A < \frac{3P}{4} &\Rightarrow u < \frac{3P\gamma - \gamma P - 1}{4} \\
    u \le v &\Rightarrow u^2 \le c^2 \\
    A \le B &\Rightarrow u \le v
\end{align*}

\subsection{Filtering admissible values of \mh{u}}\label{subsec:filter-u}

The admissible values of \( u \) must satisfy the following conditions:

\[
\begin{cases}
    u \mid c^2 \\
    u \equiv -c \pmod{\gamma} \\
    u \not\equiv -c \pmod{P} \\
    u \le v
\end{cases}
\]

Invalid values are discarded if any of the conditions fail.

\subsection{Method check for \mh{P = 2521}}\label{subsec:check-P2521}

\begin{figure}
\centering
\includegraphics[width=1 \textwidth]{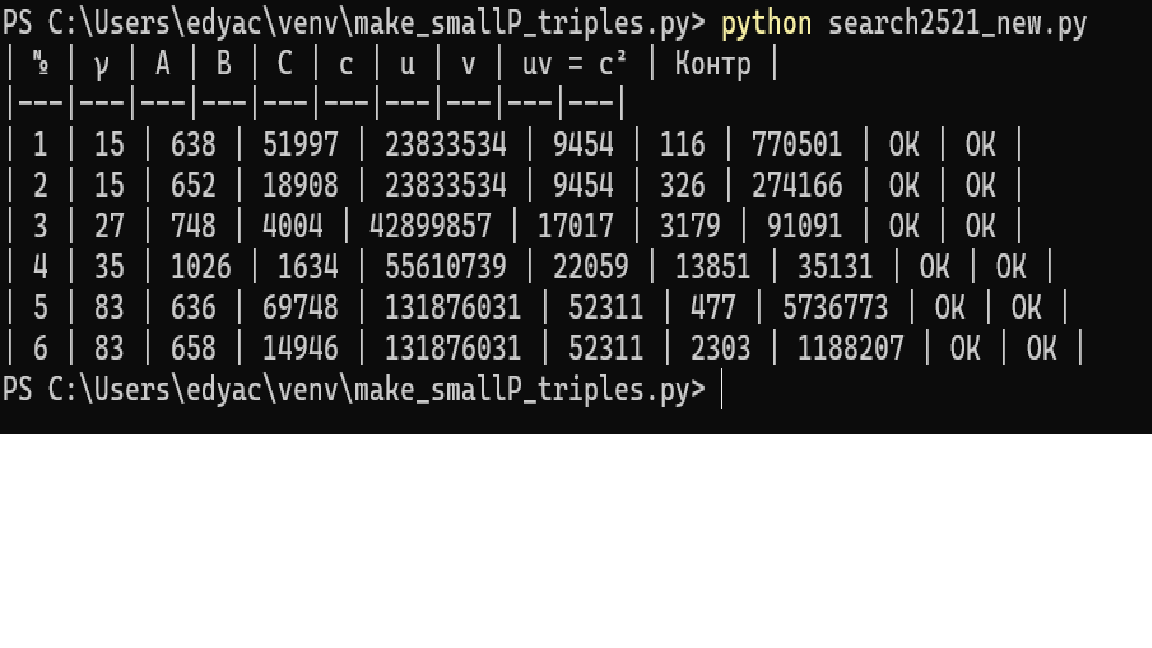}
\caption{Console output screenshot from \texttt{search2521\_new.py} for parameter $P=2521$: the table displays computed values $(A, B, C)$, coefficients $c$, $u$, the product $uv = c'$, and validation checks. All rows pass successfully (indicated by the “OK” columns).}

\label{F:ES_2521}
\end{figure}

As a result of a test run of the algorithm, six solutions were found [tab. \ref{tab:1}] and computational output [fig.\ref{F:ES_2521}]. For each of them, the parameters were computed:

\[
c = \frac{C}{P}, \quad u = \gamma A - c, \quad v = \gamma B - c
\]

Conditions $uv = c^2$, $\gamma \mid (u+c)$, $\gamma \mid (v+c)$ were verified.
\begin{table}[h!]
\centering
\caption{Enumeration results of ED1 for $P=2521$}
\begin{tabular}{|r|r|r|r|r|r|r|r|c|c|}
\hline
№ & $\gamma$ & $A$ & $B$ & $C$ & $c$ & $u$ & $v$ & $uv = c^2$ & Congr. \\
\hline
1 & 15 & 638  & 51997 & 23833534  & 9454  & 116  & 770501  & OK & OK \\
2 & 15 & 652  & 18908 & 23833534  & 9454  & 326  & 274166  & OK & OK \\
3 & 27 & 748  & 4004  & 42899857  & 17017 & 3179 & 91091   & OK & OK \\
4 & 35 & 1026 & 1634  & 55610739  & 22059 & 1851 & 35131   & OK & OK \\
5 & 83 & 636  & 69748 & 131876031 & 52311 & 477  & 5736773 & OK & OK \\
6 & 83 & 658  & 14946 & 131876031 & 52311 & 2303 & 1188207 & OK & OK \\ 
\hline
\label{tab:1}
\end{tabular}
\end{table}

In the “Congr.” column, OK means that the congruences

\[
u \equiv -c \pmod{\gamma}, \quad v \equiv -c \pmod{\gamma}.
\]

are satisfied. In all six cases, conditions $uv = c^2$ and divisibility by $\gamma$ are fully confirmed, demonstrating the correctness of the algorithm for the chosen $P$.

\subsection{Relation between different solutions}\label{subsec:c'_C2}

Consider two solutions of the same equation.
\[
\frac{4}{P}=\frac{1}{A_i}+\frac{1}{B_i}+\frac{1}{C_i},\quad
A_i\le B_i\le C_i\in\mathbb{N},\quad C_i=c_iP,\quad i=1,2,
\]
linked by the conditions;
\[
4c_i-1=\gamma_i P,\quad \gamma_i=4\nu_i-1,\quad c_i=\frac{\gamma_i P+1}{4}.
\]
Let $\nu:=\nu_2-\nu_1$ and $c_2>c'$ (that is, $\nu\ge 0$).

Introduce the notation:
\[
g_i=\gcd(A_i,B_i),\quad A_i=g_i a_i,\quad B_i=g_i b_i,\quad \gcd(a_i,b_i)=1.
\]

\begin{theorem}
The relation
\[
\frac{B_2+A_2}{A_2B_2}>1
\]
holds if and only if, under the condition
\begin{align*}
&\gcd(A_1,B_1)=g_1,\quad A_1=g_1a_1,\quad B_1=g_1b_1,\quad \gcd(a_1,b_1)=1,\\
&\gcd(A_2,B_2)=g_2,\quad A_2=g_2a_2,\quad B_2=g_2b_2,\quad \gcd(a_2,b_2)=1,\\
&\frac{b_2+a_2}{g_2a_2b_2}-\nu=\frac{b_1+a_1}{g_1a_1b_1},\quad \nu=\nu_2-\nu_1,
\end{align*}
we have $g_2a_2=1$.
\end{theorem}

\begin{proof}
We have the identity
\[
\frac{B_2+A_2}{A_2B_2}=\frac{a_2+b_2}{g_2 a_2 b_2}.
\]
($\Rightarrow$) Suppose $\frac{B_2+A_2}{A_2B_2}>1$. Then
\[
a_2+b_2>g_2 a_2 b_2
\quad\Longleftrightarrow\quad
\frac{1}{a_2}+\frac{1}{b_2}>g_2.
\]
The left-hand side $\le 1+\frac{1}{\min\{a_2,b_2\}}\le 2$, and since $g_2\in\mathbb{N}$, it follows that $g_2=1$. Furthermore,
\[
\frac{1}{a_2}+\frac{1}{b_2}>1
\quad\Longleftrightarrow\quad
\min\{a_2,b_2\}=1.
\]
If $a_2\le b_2$ we have $a_2=1$. Together with $g_2=1$, this gives $g_2a_2=1$.

($\Leftarrow$) Suppose $g_2a_2=1$, i.e. $g_2=1$, $a_2=1$. Then
\[
\frac{B_2+A_2}{A_2B_2}=\frac{1+b_2}{b_2}>1.
\]
Thus, the condition $g_2a_2=1$ is equivalent to the required inequality.
The link with $\nu$ in the formula
\[
\frac{b_2+a_2}{g_2 a_2 b_2}-\nu=\frac{b_1+a_1}{g_1 a_1 b_1},\quad \nu\in\mathbb{N},
\]
simply aligns the quantities: the left fraction with $\nu\ge 1$ and a positive right-hand side is indeed $>1$.
\end{proof}

\begin{proposition}
Let $\nu=0$ and
\[
\frac{b_2+a_2}{g_2a_2b_2}=\frac{b_1+a_1}{g_1a_1b_1}.
\]
If $\frac{B_2+A_2}{A_2B_2}>1$, then $(A_1,B_1)=(A_2,B_2)$.
\end{proposition}

\begin{proof}
From $\frac{a_2+b_2}{g_2a_2b_2}>1$ it follows that $g_2=1$, $a_2=1$. Then
\[
\frac{1+b_2}{b_2}=\frac{a_1+b_1}{g_1a_1b_1}>1
\]
gives $g_1=1$, $a_1=1$ and $b_1=b_2$. Therefore, $(A_1,B_1)=(A_2,B_2)$.
\end{proof}

\paragraph{General relations for \(\nu=0\).}
Fix \(P\), \(\gamma\), \(c\) so that \(4c-1=\gamma P\), \(\gcd(\gamma,c)=1\).
Then for any divisor \(u\mid c^2\) with \(u\equiv -c\pmod{\gamma}\) (and \(v:=c^2/u\)) we obtain the following.
\[
A=\frac{u+c}{\gamma},\qquad
B=\frac{v+c}{\gamma},\qquad
C=cP.
\]
The pair \((u,v)\) and its permutation \((v,u)\) give the same set \(\{A,B\}\).

\subsubsection*{Row pair 1–2: \(\gamma=15,\; c=9454,\; P=2521\)}
Check of invariants:
\[
4c-1=4\cdot 9454-1=37816-1=37815=\gamma P=15\cdot 2521,\quad
C=cP=9454\cdot 2521=23\,833\,534.
\]
Factorization: \(c=2\cdot 4727\), \quad \(c^2 = 2^2\cdot 4727^2\).  
Residue class for branching by \(u\):
\[
c\equiv 4\pmod{15}\ \Rightarrow\ -c\equiv 11\pmod{15}.
\]

\textbf{Row 1:} \(u=116=2^2\cdot 29\) (\(116\equiv 11\pmod{15}\)),  
\(v=\frac{c^2}{u}=\frac{89\,378\,116}{116}=770\,501\), factorization \(v=13\cdot 59\,269\).  
Check: \(u v = c^2\).  
\[
A=\frac{116+9454}{15}=638,\quad
B=\frac{770\,501+9454}{15}=51\,997,\quad
C=23\,833\,534.
\]

\textbf{Row 2:} \(u=326=2\cdot 163\) (\(326\equiv 11\pmod{15}\)),  
\(v=\frac{89\,378\,116}{326}=274\,166=2\cdot 137\,083\).  
Check: \(u v = c^2\).  
\[
A=\frac{326+9454}{15}=652,\quad
B=\frac{274\,166+9454}{15}=18\,908,\quad
C=23\,833\,534.
\]

\emph{Branching point:} same \((\gamma,c)\), one residue class \(u\equiv -c\pmod{15}\), but different \(u\) give different \(A,B\) for the same \(C\).
\subsubsection*{Row pair 5–6: \(\gamma=83,\; c=52311,\; P=2521\)}
Check of invariants:
\[
4c-1=4\cdot 52311-1=209244-1=209243=\gamma P=83\cdot 2521,\quad
C=cP=52311\cdot 2521=131\,876\,031.
\]
Factorization: \(c=3\cdot 7\cdot 47\cdot 53\), \quad \(c^2 = 3^2\cdot 7^2\cdot 47^2\cdot 53^2\).  
Residue class for branching by \(u\):
\[
c\equiv 21\pmod{83}\ \Rightarrow\ -c\equiv 62\pmod{83}.
\]

\textbf{Row 5:} \(u=477=3^2\cdot 53\) (\(477\equiv 62\pmod{83}\)),  
\(v=\frac{c^2}{u}=\frac{2\,736\,435\,?}{477}=5\,736\,773=7^2\cdot 47^2\cdot 53\).  
Check: \(u v = c^2\).  
\[
A=\frac{477+52311}{83}=636,\quad
B=\frac{5\,736\,773+52311}{83}=69\,748,\quad
C=131\,876\,031.
\]

\textbf{Row 6:} \(u=2303=7^2\cdot 47\) (\(2303\equiv 62\pmod{83}\)),  
\(v=\frac{c^2}{u}=1\,188\,207=3^2\cdot 47\cdot 53^2\).  
Check: \(u v = c^2\).  
\[
A=\frac{2303+52311}{83}=658,\quad
B=\frac{1\,188\,207+52311}{83}=14\,946,\quad
C=131\,876\,031.
\]

\emph{Branching point:} same \((\gamma,c)\), one residue class \(u\equiv -c\pmod{83}\), but different divisors \(u\) give different \(A,B\) for the same \(C\).

\paragraph{Conclusion.}
For fixed \(\gamma, c\) and \(\nu=0\) the large denominator \(C=cP\) remains unchanged, and the branching of solutions is completely determined by the set of admissible divisors.
\[
\mathcal{U}_{\gamma,c}=\{\,u\mid c^2:\;u\equiv -c\pmod{\gamma}\,\}.
\]
Each \(u\in\mathcal{U}_{\gamma,c}\) generates its own pair \((A,B)\) via the formulas above; the permutation \((u,v)\leftrightarrow (v,u)\) corresponds to the permutation \((A,B)\).

\paragraph{Case $\nu=1$: construction of the second solution and constraints.}
Let us
\[
\frac{b_2+a_2}{g_2a_2b_2}-1=\frac{b_1+a_1}{g_1a_1b_1},\qquad
\gcd(a_i,b_i)=1,\quad A_i=g_i a_i,\ B_i=g_i b_i.
\]
Then equivalently
\[
\frac{b_2+a_2}{g_2a_2b_2}
=1+\frac{b_1+a_1}{g_1a_1b_1}
=\frac{g_1a_1b_1+a_1+b_1}{g_1a_1b_1}.
\]
From this, after cross-multiplication,
\begin{equation}\label{eq:nu1-master}
(b_2+a_2)\,g_1a_1b_1
= g_2a_2b_2\,(g_1a_1b_1+a_1+b_1).
\end{equation}

\begin{proposition}[Base branch with $g_2a_2=1$]
If $g_2a_2=1$ (that is, $g_2=1$, $a_2=1$), then \eqref{eq:nu1-master} reduces to
\[
(1+b_2)\,g_1a_1b_1=b_2\,(g_1a_1b_1+a_1+b_1)\ 
\Longleftrightarrow g_1a_1b_1=b_2\,(a_1+b_1).
\]

Hence,
\[a_1+b_1 \mid g_1\quad\text{and}\quad b_2=\frac{g_1a_1b_1}{a_1+b_1}.
\]
In terms of $(A,B)$ this gives
\[
A_2=1,\qquad
B_2=\frac{g_1a_1b_1}{a_1+b_1},
\]
and $B_2\in\mathbb{N}$ if and only if $a_1+b_1\mid g_1$.
\end{proposition}

\begin{remark}
The general form of \eqref{eq:nu1-master} can be solved for $b_2$:
\[
b_2=\frac{a_2\,g_1a_1b_1}{\,g_2a_2\,(g_1a_1b_1+a_1+b_1)-g_1a_1b_1\,}.
\]
The denominator must be a positive divisor of the numerator. The branch $g_2a_2=1$ is the lowest in $A_2$ and the only one consistent with the criterion
\(
\frac{A_2+B_2}{A_2B_2}>1\iff g_2a_2=1
\)
(the theorem above). In other words, for $\nu=1$ the 'second' solution, if it exists, necessarily has $A_2=1$.
\end{remark}

\begin{proposition}[Compatibility with the original equation]
For any pair $(A,B)$, the value $c$ is recovered from the identity
\[
\frac{A+B}{AB}=\frac{4}{P}-\frac{1}{cP}
\iff
c=\frac{AB}{\,4AB-P(A+B)\,}.
\]
For $A_2=1$, $B_2=b_2$ we obtain
\[
c_2=\frac{b_2}{\,4b_2-P(1+b_2)\,}
=\frac{b_2}{\,b_2(4-P)-P\,}.
\]
The requirement $c_2\in\mathbb{N}$, $c_2>0$ implies the condition
\[
4b_2-P(1+b_2)>0
\iff
b_2(4-P)>P.
\]
This yields a strict restriction on the prime $P$:
\[
P\le 3.
\]
For $P\ge 5$, the inequality is impossible for any $b_2\in\mathbb{N}$, i.e., the “branch” $\nu=1$ does not give an admissible second solution within the original equation (with $C=cP$).
\end{proposition}

\paragraph{Procedure for $\nu=1$ (generalization and filters).}
\begin{enumerate}
  \item \textbf{Normalization of the first solution:} Find $g_1=\gcd(A_1,B_1)$, $A_1=g_1a_1$, $B-1=g_1b_1$, $\gcd(a_1,b_1)=1$.
  \item \textbf{Divisibility check:} If $a_1+b_1\nmid g_1$, the branch $\nu=1$ with $g_2a_2=1$ is impossible.
  \item \textbf{Candidate reconstruction:} If $a_1+b_1\mid g_1$, set
  \[
  A_2=1,\qquad B_2=\frac{g_1a_1b_1}{a_1+b_1}.
  \]
  \item \textbf{Compatibility with $P$:} Compute
  \[
  c_2=\frac{A_2B_2}{\,4A_2B_2-P(A_2+B_2)\,}=\frac{b_2}{\,4b_2-P(1+b_2)\,}.
  \]
  Require $c_2\in\mathbb{N}$ and $c_2>0$. This is possible only for $P\in\{2,3\}$.
  \item \textbf{Ordering $A_2\le B_2\le C_2$:} 
  \[
  A_2=1\le B_2,\qquad B_2\le C_2=c_2P
  \]
  holds automatically when $c_2>0$; for the strict inequality $A_2<B_2$ one needs $b_2\ge 2$.
\end{enumerate}
\subsection{Case \mh{B \equiv 0 \pmod P; A,C \not\equiv 0 \pmod P}}\label{subsec:case-B-mult-P}

Similarly to the analysis in~\S\ref{sec:param-Bnot0}, the substitution \(B = bP\) with

\[
4b - 1 = \gamma P,\qquad b = \frac{\gamma P + 1}{4},\quad \gcd(\gamma,b)=1
\]

leads to the identity (in parametrization form);

\[
  (\gamma A - b)\,(\gamma C - b) = b^2.
\]

Parameterization is carried out via a pair of divisors \(u,v\) of the number \(b^2\):

\[
  uv = b^2,\qquad u \equiv v \equiv -b \pmod{\gamma},
\]

with the formulas

\[
  A = \frac{u+b}{\gamma},\qquad B = bP,\qquad C = \frac{v+b}{\gamma}.
\]

The modulus‑\(P\) filter is now imposed on \(C\):

\[
  C \not\equiv 0 \pmod P \;\Longleftrightarrow\; v \not\equiv -b \pmod P.
\]

The enumeration algorithm, normalization, and duplicate elimination are completely analogous to \mbox{Case~1}, with the substitutions \(c\to b\) and \(B\leftrightarrow C\) in the notation.

Additionally, within the ordering \(A \le B \le C\) one requires

\[
  B < C \;\Longleftrightarrow\; bP < \frac{v+b}{\gamma},
\]

which is equivalent to \(\gamma bP < v+b\). When enumerating the admissible \(u,v\), this condition serves to cut off configurations where the second and third denominators coincide or violate the order.

\subsection{Absence of ordered solutions for \mh{ B \equiv 0 \pmod P, A,C \not\equiv 0 \pmod P }}\label{subsec:no-sol-B-mult-P}

In the notation of the previous paragraph (\(B=bP\), \(4b-1=\gamma P\), \(b=(\gamma P+1)/4\)), the ordering condition \(A \le B \le C\) requires

\[
  B \le C \;\Longleftrightarrow\; bP \le \frac{v+b}{\gamma} \le \frac{b^2+b}{\gamma}.
\]

Hence,

\[
  bP \le \frac{b^2+b}{\gamma}
  \;\Longleftrightarrow\;
  \gamma b P \le b^2 + b
  \;\Longleftrightarrow\;
  \gamma P \le b + 1.
\]

Substituting \(b=\frac{\gamma P+1}{4}\), we obtain the following.

\[
  \gamma P \le \frac{\gamma P+1}{4} + 1 = \frac{\gamma P+5}{4}
  \;\Longleftrightarrow\;
  4\gamma P \le \gamma P + 5
  \;\Longleftrightarrow\;
  3\gamma P \le 5.
\]

For any odd prime \(P\) and \(\gamma \ge 1\) this is impossible. Therefore, for odd \(P\), the standard ordering \(A \le B \le C\) is incompatible with the subcase \(B \equiv 0 \pmod P\), \(A,C \not\equiv 0 \pmod P\).

\begin{remark}
It follows that in all solutions where exactly one denominator is divisible by \(P\), this denominator must be the largest one, i.e., \(C=cP\) (see “Case~1”). The subcase with \(B \equiv 0 \pmod P\) under ordering reduces to “Case~1” by a simple permutation.
\end{remark}

\section{Case of two multiples of \mh{P: B = bP, C = cP, A \not\equiv 0 \pmod P}}\label{sec:two-multiples}

Consider
\[
  \frac{4}{P} = \frac{1}{A} + \frac{1}{bP} + \frac{1}{cP},
  \qquad A \le bP \le cP, \quad b,c \in \mathbb{N}.
\]
Multiplying by \(A\,b\,c\,P\), we obtain the following.
\[
  4A b c = P b c + A(b + c)
  \quad\Longleftrightarrow\quad
  A\,(4bc - b - c) = P b c.
\]
Denote
\[
  t := 4bc - b - c,
\]
then
\[
  A = \frac{P b c}{t}.
\]
Since \(t > bc\) for \(b,c \ge 1\), for \(A\) to be an integer, it is necessary that \(P \mid t\). Set
\[
  t = P\,\delta, \qquad \delta \in \mathbb{N}.
\]

\begin{lemma}\label{lem:bc_divisibility}
If \(4bc - b - c = P\delta\), then \(\delta \mid bc\), and
\[
  A \;=\; \frac{bc}{\delta}.
\]
\end{lemma}
\begin{proof}
From \(A(4bc - b - c) = P\,bc\) and \(4bc - b - c = P\delta\) we get \(A\,P\delta = P\,bc\).
By cancelation \(P\), we have \(A\,\delta = bc\), so \(\delta \mid bc\) and \(A = bc/\delta\).
\end{proof}
\begin{lemma}[Factorization via \(\delta=\alpha d'^2\)]\label{lem:delta-split}
Let \(b,c\in\mathbb{N}\), \(g:=\gcd(b,c)\), \(b=gb'\), \(c=gc'\), \(\gcd(b',c')=1\),
and \(4bc-b-c=P\delta\).
Represent \(\delta\) in the form \(\delta=\alpha d'^2\), where \(\alpha\) is square‑free, \(d'\in\mathbb{N}\),
and set \(g:=\alpha d'\).
Then the following conditions are equivalent:
\[
\text{\rm(i)}\quad \delta \;=\; \alpha\,d'^2,\qquad
\text{\rm(ii)}\quad P d' \;=\; 4\alpha d'\,b' c' - (b'+c').
\]
In this case the identity
\[
  (4\alpha d' b' - 1)(4\alpha d' c' - 1) \;=\; 4\alpha P d'^2 + 1
\]
is valid, and in particular \(d'\mid (b'+c')\).
\end{lemma}

\begin{proof}
We write
\[
  4bc - b - c
  \;=\;
  4g^2\,b' c' - g(b'+c')
  \;=\;
  g\bigl(4g\,b' c' - (b'+c')\bigr),
\]
and for \(g=\alpha d'\) we have
\[
  P\,\delta \;=\; \alpha\,d'\,\bigl(4\alpha\,d'\,b' c' - (b'+c')\bigr).
\]
(i) \(\Rightarrow\) (ii): if \(\delta = \alpha\,{d'}^{2}\), cancel \(\alpha d'\) to obtain
\(P d' = 4\alpha d' b' c' - (b'+c')\).

(ii) \(\Rightarrow\) (i): multiplying by \(\alpha d'\) and dividing by \(P\) gives \(\delta = \alpha d'^2\).

Next, add \(\tfrac{1}{4\alpha d'}\) to both sides of (ii):
\[
  P d' + \frac{1}{4\alpha d'} \;=\;
  \Bigl(2\sqrt{\alpha d'}\,b'-\tfrac{1}{2\sqrt{\alpha d'}}\Bigr)
  \Bigl(2\sqrt{\alpha d'}\,c'-\tfrac{1}{2\sqrt{\alpha d'}}\Bigr).
\]
Multiplying by \(4\alpha d'\) yields
\[
  (4\alpha d' b' - 1)(4\alpha d' c' - 1) \;=\; 4\alpha P d'^2 + 1.
\]
From (ii) it also follows that \(b'+c' = d'(4\alpha b' c' - P)\), hence \(d' \mid (b'+c')\), and from
\(b'+c' \mid P d'\) we get \(\frac{b'+c'}{\gcd(b'+c',P)}\mid d'\).
\end{proof}

\begin{theorem}\label{thm:ED2-param}
Let \(P\) be prime and \(\delta=\frac{4bc-b-c}{P}\).
Represent \(\delta=\alpha d'^2\) with square‑free \(\alpha\) and \(d'\in\mathbb{N}\), and set \(g:=\alpha d'\).
Then all solutions with \(B=bP,\ C=cP\) are described by factorizations
\[
  X\,Y=4\alpha Pd'^2+1,\quad X\equiv Y\equiv -1\pmod{4\alpha d'},
\]
for which \(b=\frac{X+1}{4\alpha d'}\), \(c=\frac{Y+1}{4\alpha d'}\), \(\gcd(b',c')=1\) and
\(\frac{b'+c'}{\gcd(b'+c',P)}\mid d'\). Moreover,
\[
  b = g\,b',\quad c = g\,c',\quad A=\frac{bc}{\delta}=\alpha \,b'c',\quad
  B=bP,\ C=cP,\ d'\mid(b'+c'), \quad
\frac{b'+c'}{d'} \equiv 3 \pmod{4}.
\]
\end{theorem}

\subsection{Enumeration algorithm}
For a fixed prime \(P\):
\begin{enumerate}
  \item Choose \(\alpha,d'\) (e.g., by enumeration); when using the factorization approach — focus on factorizations \(N=4\alpha P{d'}^2+1\) and the filters of Lemma~\ref{lem:delta-split}.
  \item Form \(N = 4\alpha\,P\,d'^2 + 1\) and factor it into pairs \(X \cdot Y = N\) with \(X \equiv Y \equiv -1 \pmod{4\alpha\,d'}\).
  \item Recover \(b', c'\), check \(\gcd(b',c') = 1\) and the ordering condition \(A \le bP \le cP\).
  \item Compute \(b, c, A, B, C\) as in the theorem.
  \item Eliminate duplicates arising from swapping \(X \leftrightarrow Y\).
\end{enumerate}
\emph{Warning.} The factorization step for \(N\) is the main time‑consuming part; in practice ECM, Pollard–Rho, etc., are used.

\subsection{Example (enhanced verification for prime \mh{P=4P'+1} and condition \mh{B<C})}
\label{subsec:two-mult-example-strong-B-less-C}

Below, everything is built directly via the parameters of the main lemma $\alpha,d',b',c'$; choose $b'<c'$ so that $b<c$ and, consequently, $B=bP<C=cP$.

\paragraph{Template.}
Let $P\equiv 1\pmod{4}$ be the prime, $\alpha,d'\in\mathbb{N}$, and suppose that there exist $b',c'\in\mathbb{N}$ such that
\[
(4 \alpha d' b'-1)\,(4 \alpha d' c'-1)=4 \alpha P {d'}^{2}+1,\qquad b'<c'.
\]
Set
\[
g:=\alpha d',\quad b:=g\,b',\quad c:=g\,c',\quad \delta:=\alpha \,{d'}^{2}.
\]
Then
\[
4bc-b-c=\delta\,P,\qquad
A=\frac{bc}{\delta}=\alpha \,b' c',\quad B=bP,\quad C=cP,
\]
and the decomposition
\[
\frac{4}{P}=\frac{1}{A}+\frac{1}{B}+\frac{1}{C},\qquad B<C
\]
holds.

\paragraph{Example 1: $P=29=4\cdot 7+1$, $\alpha =1$, $d'=2$.}
Here
\[
N=4 \alpha P {d'}^{2}+1=4\cdot 1\cdot 29\cdot 4+1=465=15\cdot 31.
\]
Take the ordered parameters.
\[
4 \alpha d' b'-1=15,\quad 4 \alpha d' c'-1=31
\quad\Longrightarrow\quad
b'=\frac{15+1}{8}=2,\ \ c'=\frac{31+1}{8}=4,
\]
so that $b'<c'$. Then
\[
g=\alpha d'=2,\quad b=g b'=4,\quad c=g c'=8,\quad \delta=\alpha{d'}^{2}=4.
\]
Check:
\[
4bc-b-c=4\cdot 4\cdot 8-4-8=128-12=116=\delta P=4\cdot 29.
\]
Denominators:
\[
A=\frac{bc}{\delta}=\frac{32}{4}=8,\quad
B=bP=4\cdot 29=116,\quad
C=cP=8\cdot 29=232,
\]
in fact $B<C$, as well as
\[
\frac{4}{29}=\frac{1}{8}+\frac{1}{116}+\frac{1}{232}.
\]

\paragraph{Example 2: $P=53=4\cdot 13+1$, $\alpha=1$, $d'=3$.}
Here
\[
N=4 \alpha P {d'}^{2}+1=4\cdot 1\cdot 53\cdot 9+1=1909=23\cdot 83.
\]
Parameters
\[
4 \alpha d' b'-1=23,\quad 4 \alpha d' c'-1=83
\quad\Longrightarrow\quad
b'=\frac{23+1}{12}=2,\ \ c'=\frac{83+1}{12}=7,
\]
and $b'<c'$. Then
\[
g=\alpha d'=3,\quad b=g b'=6,\quad c=g c'=21,\quad \delta=\alpha{d'}^{2}=9.
\]
Check:
\[
4bc-b-c=4\cdot 6\cdot 21-6-21=504-27=477=\delta P=9\cdot 53.
\]
Denominators:
\[
A=\frac{bc}{\delta}=\frac{126}{9}=14,\quad
B=bP=6\cdot 53=318,\quad
C=cP=21\cdot 53=1113,
\]
and $B<C$, as well as
\[
\frac{4}{53}=\frac{1}{14}+\frac{1}{318}+\frac{1}{1113}.
\]

\subsection{Verification of the new algorithm for \mh{P = 2521,3529} (case \mh{B = bP, C = cP})}\label{subsec:two-mult-check-P2521}

In a test run of the algorithm, solutions were found.  
For each, the parameters were computed:
\[
  \delta = \frac{4bc - b - c}{P},\quad
  X = 4\alpha d' b_{1} - 1,\quad
  Y = 4\alpha d' c_{1} - 1
\]
and the conditions were checked [tab.\ref{tab:2}]:
\[
  X\,Y = 4\alpha P\,{d'}^{2} + 1,\quad \delta \mid bc,
\]
as well as the order \(A < B \le C\) and the main equality.
\[
  \frac{4}{P} = \frac{1}{A} + \frac{1}{B} + \frac{1}{C}.
\]


\begin{table}[h!]
\centering
\caption{Results for $P = 2521$ and $P = 3529$ with verification by the lemma $X \cdot Y = 4 \alpha P (d')^2 + 1$}
\begin{tabular}{|c|c|c|c|c|c|c|c|c|c|c|c|c|c|c|c|}
\hline
\# & $\alpha$ & $b'$ & $c'$ & $g$ & $b$ & $c$ & $\delta$ & $X$ & $Y$ & $N$ & $A$ & $B$ & $C$ & $d'$ & OK \\ \hline
\multicolumn{16}{|c|}{$P = 2521$} \\ \hline
1 & 1 & 4 & 161 & 3 & 12 & 483 & 9 & 47 & 1931 & 90757 & 644 & 30252 & 1217643 & 3 & \checkmark \\
2 & 2 & 2 & 159 & 14 & 28 & 2226 & 98 & 111 & 8903 & 988233 & 636 & 70588 & 5611746 & 7 & \checkmark \\
3 & 11 & 2 & 29 & 11 & 22 & 319 & 11 & 87 & 1275 & 110925 & 638 & 55462 & 804199 & 1 & \checkmark \\ \hline
\multicolumn{16}{|c|}{$P = 3529$} \\ \hline
1 & 1 & 5 & 186 & 1 & 5 & 186 & 1 & 19 & 743 & 14117 & 930 & 17645 & 656394 & 1 & \checkmark \\
2 & 1 & 3 & 307 & 2 & 6 & 614 & 4 & 23 & 2455 & 56465 & 921 & 21174 & 2166806 & 2 & \checkmark \\
3 & 1 & 3 & 296 & 13 & 39 & 3848 & 169 & 155 & 15391 & 2385605 & 888 & 137631 & 13579592 & 13 & \checkmark \\
4 & 2 & 4 & 111 & 10 & 40 & 1110 & 50 & 159 & 4439 & 705801 & 888 & 141160 & 3917190 & 5 & \checkmark \\
5 & 5 & 1 & 181 & 10 & 10 & 1810 & 20 & 39 & 7239 & 282321 & 905 & 35290 & 6387490 & 2 & \checkmark \\
6 & 13 & 4 & 17 & 39 & 156 & 663 & 117 & 623 & 2651 & 1651573 & 884 & 550524 & 2339727 & 3 & \checkmark \\
7 & 17 & 2 & 26 & 68 & 136 & 1768 & 272 & 543 & 7071 & 3839553 & 884 & 479944 & 6239272 & 4 & \checkmark \\
8 & 26 & 1 & 34 & 130 & 130 & 4420 & 650 & 519 & 17679 & 9175401 & 884 & 458770 & 15598180 & 5 & \checkmark \\ \hline
\label{tab:2}
\end{tabular}
\end{table}

\noindent\textbf{Verification comments:}
\begin{itemize}
    \item For each row, the check of Theorem \ref{thm:ED2-param} was performed:
    \[
        X \cdot Y = 4 \cdot \alpha \cdot P \cdot (d')^2 + 1
    \]
    The equality was verified in integers - in all cases it is true (OK = \checkmark).
    \item In the rows for $P = 2521$ and $P = 3529$ there is full agreement with the constructive formulas of the algorithm:
    \[
        g = \alpha \cdot d',\quad b = g \cdot b',\quad c = g \cdot c',\quad \delta = \alpha (d')^2
    \]
    as well as
    \[
        A = \alpha b' c',\quad B = b \cdot P,\quad C = c \cdot P.
    \]
    \item Example check ($P = 2521$, row 1):
    \[
        d' = 3,\quad N = 4 \cdot 1 \cdot 2521 \cdot 3^2 + 1 = 90757,
    \]
    \[
        X \cdot Y = 47 \cdot 1931 = 90757
    \]
    — the equality holds.
    \item Example check ($P = 3529$, row 6):
    \[
        d' = 3,\quad N = 4 \cdot 13 \cdot 3529 \cdot 3^2 + 1 = 1651573,
    \]
    \[
        X \cdot Y = 623 \cdot 2651 = 1651573
    \]
    — the equality holds.
\end{itemize}
In all solutions, the following conditions hold simultaneously:
\(X \cdot Y = 4\alpha P d'^{2} + 1\) (parametrization identity; see Theorem~\ref{thm:ED2-param}),
\(\delta \mid bc\) (by Lemma~\ref{lem:bc_divisibility}),
as well as the ordering \(A < B \le C\) and the main equality
\(\frac{4}{P} = \frac{1}{A} + \frac{1}{B} + \frac{1}{C}\).
The coincidence of all computed values with the theoretical ones
confirms the correctness of the new algorithm for the chosen $P$.

\section{Impossibility of other multiplicity configurations}\label{sec:impossible-multiplicity}

\begin{lemma}\label{lem:at-least-one-mult}
In any solution, at least one of the denominators $(A,B,C)$ is divisible by~\(P\).
\end{lemma}

\begin{proof}
Multiplying the original equality by \(4PABC\), we get
\[
 4ABC = P(AB + AC + BC)
\]
Since $P$ is prime, $P$ divides at least one of $A,B,C$.
\end{proof}

\begin{lemma}\label{lem:three-mults-impossible}
Three simultaneous multiplicities are impossible: \(A = aP\), \(B = bP\), \(C = cP\).
\end{lemma}
\begin{proof}
After canceling $P$, we have
\[
  4 = \frac{1}{a} + \frac{1}{b} + \frac{1}{c} \le 3,
\]
and since $a,b,c \in \mathbb{N}$ we obtain a contradiction.
\end{proof}

\begin{lemma}[Single multiple denominator]\label{lem:one-mult-is-C}
Let in the decomposition
\[
  \frac{4}{P} = \frac{1}{A} + \frac{1}{B} + \frac{1}{C}
\]
exactly one denominator be divisible by \(P\). Then (after reordering) it is \(C = cP\).
\end{lemma}

\begin{proof}
By symmetry in $(A,B,C)$ we may assume without loss of generality that \(A \le B \le C\).

The case \(B = bP\) with \(P \nmid A\) and \(P \nmid C\) has already been excluded (see~\S\ref{subsec:no-sol-B-mult-P}): it leads to the impossible inequality \(3\gamma P \le 5\).

It remains to consider \(A = aP\) with \(P \nmid B\) and \(P \nmid C\). Then
\[
  \frac{1}{B} + \frac{1}{C} 
  = \frac{4}{P} - \frac{1}{aP} 
  = \frac{4a - 1}{aP}.
\]
Since \(P\) does not divide \(B\) and \(C\), the left-hand side in lowest terms has a denominator not divisible by \(P\), hence in the right-hand side the factor \(P\) in the denominator must cancel. Thus \(P \mid (4a - 1)\), i.e., \(4a - 1 = kP\) with \(k \in \mathbb{N}\).

From the ordering \(A = aP \le B \le C\) it follows that
\[
  \frac{1}{B} \le \frac{1}{aP},\quad \frac{1}{C} \le \frac{1}{aP},
\]
and
\[
  \frac{k}{a} = \frac{1}{B} + \frac{1}{C} \le \frac{2}{aP}
  \quad\Longrightarrow\quad k \le \frac{2}{P}.
\]
For \(P \ge 5\) we have \(2/P < 1\), hence \(k = 0\), which is impossible. Therefore, the case \(A = aP\) is excluded, and with one multiple it must be \(C\).
\end{proof}

\begin{remark}\label{rem:one-multiple-order}
The key point is the right to order the denominators \(A \le B \le C\) without loss of generality — this follows from the symmetry of the original equation and justifies the use of estimates like \(1/B \le 1/A\).
\end{remark}

\begin{lemma}\label{lem:two-multiples-BC}
If exactly two denominators are divisible by~\(P\), then they must be \(B = bP\) and \(C = cP\); the configurations \(A = aP,\,B = bP,\,P\nmid C\) or \(A = aP,\,C = cP,\,P\nmid B\) are impossible.
\end{lemma}

\begin{proof}
Multiplying the original equality by \(P\), we get
\[
  4 = \frac{P}{A} + \frac{P}{B} + \frac{P}{C}.
\]
If \(A\) and \(B\) are divisible by \(P\) but \(C\) is not, then the first two terms are integers and the third not, so the sum cannot be an integer. The same holds for the pair \((A,C)\). Therefore, the two multiples can only be \(B\) and \(C\).
\end{proof}
\section{Algorithm convergence and conditional completeness of coverage}
\label{sec:conv-completeness}

\subsection{Logarithmic density of parameters}
\label{subsec:log-density}

\begin{theorem}[see also Theorem~\ref{th:log_convergence}]
\label{th:log_density}
Let $\Lambda \subset \mathbb{Z}^k$ be a sublattice of rank~$k$ with index~$M$ independent of the prime~$P$. Consider the parametric box
\[
  \mathcal{B}_k(T) = \{ u \in \mathbb{Z}^k \mid 1 \le u_i \le T \},
\]
and the set of admissible parameters
\[
  \mathcal{G}_P(T) = \{ u \in \mathcal{B}_k(T) : u \equiv u_0(P) \pmod{\Lambda} \}.
\]
Then
\[
  |\mathcal{G}_P(T)| = \frac{T^k}{M} + O_k\!\left(T^{k-1}\right),
\]
and this asymptotic remains valid for $T = (\log P)^A$ for any fixed $A > 0$.
\end{theorem}

\begin{proof}
Since $\Lambda$ is a lattice of index~$M$, the number of residue classes modulo~$\Lambda$ in $\mathcal{B}_k(T)$ is equal to $\frac{T^k}{M}$ with an error term $O_k(T^{k-1})$. The residue $u_0(P)$ specifies exactly one residue class, so for $T = (\log P)^A$ we obtain logarithmic growth of the size of the set~$\mathcal{G}_P(T)$.
\end{proof}

\subsection{Logarithmic convergence of the algorithm}
\label{subsec:log-convergence}

\begin{theorem}
\label{th:log_convergence}
Let the admissible parameters $u$ lie in the affine class $u \equiv u_0(P) \pmod{\Lambda}$, where $\Lambda \subset \mathbb{Z}^k$ is a sublattice of fixed index~$M$. Consider an algorithm that enumerates all $u \in \mathcal{B}_k\big((\log P)^A\big)$ and tests admissibility. Then:
\begin{enumerate}
  \item There exists a constant $c>0$, depending on~$\Lambda$, such that
    \[
      \big|\mathcal{G}_P\big((\log P)^A\big)\big| \ge c \cdot \frac{(\log P)^{Ak}}{M}.
    \]
  \item The average number of iterations before finding an admissible parameter is bounded above by a constant independent of~$P$.
  \item A full search is guaranteed to find a solution in
    \[
      O\big((\log P)^{Ak}\big)
    \]
    steps.
\end{enumerate}
\end{theorem}

\begin{proof}
This follows from Theorem~\ref{th:log_density}, since in the enumerated box the admissible parameters occupy a positive proportion. Due to the uniformity of the distribution and the constancy of the lattice index, the average time to success is bounded, and a full search covers $O((\log P)^{Ak})$ points.
\end{proof}

\begin{corollary}
For any fixed $A > 0$ and search radius
\[
T = (\log P)^A
\]
the enumeration algorithm of \S9.2 will find an ED2 solution in $O\big((\log P)^{3A}\big)$ steps on average.
\end{corollary}

\begin{remark}[On the limits of applicability of Theorems 9.1–9.2]
These theorems estimate the number of points on an affine lattice in parametric boxes, but by themselves
do not guarantee the existence of solutions to the nonlinear identity \((4b-1)(4c-1)=4P\delta+1\).
To connect the “geometry of enumeration” with the existence of a triple \((\delta,b,c)\) an external input is required —
averaging over \(\delta\) (BV-type estimates for \(S(\delta)\)) and/or a construction via the parametrization \((t,k)\) (§9.10),
which we use later.
\end{remark}

\subsection{Affine class of parameters}
\label{subsec:affine-class}

\begin{definition}
\label{def:affine-class}
Let $\Lambda \subset \mathbb{Z}^k$ be a full sublattice.  
An \emph{affine class of parameters} is the set
\[
  \mathcal{G}_P = u_0(P) + \Lambda
    = \{ u \in \mathbb{Z}^k : u \equiv u_0(P) \pmod{\Lambda} \},
\]
describing compatible sets of parameters satisfying the local conditions of the algorithm
(integrality, co-primality, modular constraints).
\end{definition}

\subsection{Parametric boxes}
\label{subsec:param-boxes}

\subsubsection{Type I box}\label{subsec:boxI}
\begin{definition}[Type I box]\label{def:boxI}
For fixed $A > 0$ and prime~$P$:
\[
  \mathcal{B}^{(I)}_P =
  \left\{ u \in \mathbb{Z}^3 :
  \begin{array}{l}
    1 \le u_i \le (\log P)^A, \\
    u \equiv u_0(P) \pmod{\Lambda_1}
  \end{array}
  \right\},
\]
where $\Lambda_1 \subset \mathbb{Z}^3$ is a sublattice of index~$M_1$ defining the modular constraints.
\end{definition}
\subsubsection{Type~II box}
\label{subsubsec:box-II}

\begin{definition}[Type~II box]
\label{def:box_II}
Let $k \ge 2$, $A > 0$, $T > 1$, and $\Lambda_2 \subset \mathbb{Z}^k$ be a sublattice of fixed index.  
A Type~II box is defined as
\[
  \mathcal{B}^{(II)}_P(T) =
  \left\{ u \in \mathbb{Z}^k :
  \begin{array}{l}
    u \equiv u_0(P) \pmod{\Lambda_2}, \\
    \rho_W(u) \le (\log P)^A, \\
    |F(u)| \le \Delta(T)
  \end{array}
  \right\},
\]
where $W \succ 0$ is a weight matrix for the norm $\rho_W(u) = \sqrt{\langle W u, u \rangle}$,  
$F \colon \mathbb{Z}^k \to \mathbb{Z}$ is a fixed (usually quadratic) form,  
$\Delta(T)$ is the window thickness.
\end{definition}

\begin{remark}
A Type~II box is used to localize parameters near nonlinear dependencies
(quadratic, bilinear, etc.), complementing the Type~I box, which covers linear classes.
\end{remark}

\begin{example}[Thickening of a quadratic surface]
For $k=3$, $u = (\delta ,b,c)$ and
\[
  F(\delta ,b,c) = (4b-1)(4c-1) - 4P \delta  - 1,
\]
we obtain
\[
  \mathcal{B}^{\mathrm{ED2}}_{\mathrm{II}}(T) =
  \left\{
  \begin{array}{l}
    (d,b,c) \equiv u_0(P) \pmod{\Lambda_2}, \\
    \rho_W(d,b,c) \le (\log P)^A, \\
    |F(d,b,c)| \le \Delta(T)
  \end{array}
  \right\}.
\]
\end{example}

\subsubsection{Radial regions}
\label{subsubsec:radial-box}

\begin{definition}[Radial box]
\label{def:radial-box}
For $W \succ 0$ and $R(T) = (\log P)^A$:
\[
  \mathcal{B}^{\mathrm{rad}}(T) =
  \left\{ u \in u_0(P) + \Lambda_2 :
    \rho_W(u) \le R(T),\;
    \Theta(u) \in \mathcal{A}
  \right\},
\]
where $\Theta(u)$ is the angular projection, restricted to a set $\mathcal{A} \subset \mathbb{S}^{k-1}$.
\end{definition}

\begin{remark}
A radial cut-off point is used as a geometric filter and can be combined
with the condition $|F(u)| \le \Delta(T)$ to localize the search.
\end{remark}

\begin{example}[Radial localization]
For $(y,c,u,v)$ with $uv = c^2$ define
\[
  \rho^2(u,v) = (\log u - \log v)^2, \quad
  \rho(u,v) \le (\log P)^{-C},
\]
which singles out nearly diagonal pairs $u \approx v$ while preserving modular conditions.
\end{example}

\begin{remark}[Status]
Type~II boxes and radial windows are not used in the proofs of Sections~9.1–9.10 and serve as groundwork for numerical experiments and future estimates.
\end{remark}

\subsection{Construction of a search algorithm on an affine lattice}
\label{subsec:lattice-search}

\paragraph{General setting.}
Let the parameter vector $u \in \mathbb{Z}^k$ satisfy a system of linear congruences
\[
  M u \equiv r \pmod{m},
\]
where $M$ is an integer matrix, $r$ is a vector of residues, and $m$ is a modulus or a set of moduli.  
The set of integer solutions to this system is described as
\[
  u_0 + \Lambda, \quad
  \Lambda = \{\, v \in \mathbb{Z}^k \mid M v \equiv 0 \pmod{m} \,\},
\]
where $u_0$ is a particular solution and $\Lambda$ is a subgroup of $\mathbb{Z}^k$. Such a set is called
an \emph{affine lattice} in~$\mathbb{Z}^k$. Any condition of the form $M u \equiv r \pmod{m}$ defines exactly
such a shift of a sublattice, and for it the results on lattice point density, minimal
vectors, etc., are applicable.

This property is \emph{central} in our convergence analysis: without the structure of an affine lattice
it is impossible to correctly apply Theorems~\ref{th:log_density}–\ref{th:log_convergence}.
\paragraph{Method.}
For fixed~$P$ and a subcase chosen from the classification, condition ED2 has the following form.

\begin{lemma}
Let $r$ be a fixed class. Then the ED2 condition is equivalent to a single linear congruence
\[
  \langle a_r, t \rangle \equiv b_r \pmod{m},
\]
where $a_r,\,b_r$ depends only on $r$ and not on $t$. The set of integer solutions for fixed~$r$ forms an affine subspace $\mathcal{A}_r \subset \mathbb{Z}^k$, and the overall set of solutions is a disjoint union of such subspaces.
\end{lemma}

For fixed~$P$ and a chosen subcase of the classification, the system of conditions on the parameters
$(\delta,b,c)$ (see Section~\ref{sec:two-multiples}) reduces to a system of linear congruences
\[
  M u \equiv r \pmod{m},
\]
where $u \in \mathbb{Z}^k$ is the parameter vector, $M$ and~$r$ are given by formulas~(...), and~$m$ is the common modulus. The set of integer solutions of this system forms an affine class.
\[
  u_0(P) + \Lambda_j,
\]
where $\Lambda_j \subset \mathbb{Z}^k$ is the subgroup defined by the homogeneous part of the system %
\footnote{In a more general setting with other constraints, the set of solutions may be an algebraic variety rather than a lattice.}.

Thus, the set of admissible parameters has the structure of an affine lattice in~$\mathbb{Z}^k$,
which ensures the applicability of Theorems~\ref{th:log_density}–\ref{th:log_convergence}
on point density.

\medskip
\noindent\textit{Note.}
The class $u_0(P) + \Lambda_j$ constructed above coincides exactly with
the affine lattice from Definition~\ref{def:affine-class},
on whose properties Theorems~\ref{th:log_density}–\ref{th:log_convergence} rely.
It is precisely to fix the corresponding component that the index~$j$ was introduced,
so these results close the proof scheme
begun in~\S\ref{sec:two-multiples} and prepare the ground
for the direct application of these theorems.
\medskip

\begin{remark}
The general method proposed above is also applicable to the case of \emph{box~I} (see~\S\,9.6),
in which the parameters $(\delta ,b,c)$ satisfy additional constraints
defining the corresponding affine class in $\mathbb{Z}^3$.
\end{remark}

\begin{enumerate}
  \item \textbf{Form} the lattice~$\Lambda_j$ and the shift~$u_0(P)$ from the modular conditions.
  \item Enumerate $u$ in the box~$\mathcal{B}^{(j)}_P$ (see Definitions \ref{def:boxI} and \ref{def:box_II}).
  \item Select parameters passing the integrality and coprimality conditions.
  \item Apply the subcase filters (divisibility, ordering of denominators, etc.).
  \item Reconstruct $(A,B,C)$ and verify the original equation.
\end{enumerate}
\subsection*{Affine lattice for ED2: explicit system based on §7.2}

We start from the structure of §7.2: let $\alpha$ be square‑free, $d'\in\mathbb{N}$, set
\[
g:=\alpha  d',\qquad \delta:=\alpha (d')^{2},\qquad b=g\,b',\ \ c=g\,c',\ \ \gcd(b',c')=1,
\]
and suppose the identity
\[
(4g b'-1)(4g c'-1)=4P\,\delta+1
\]
holds. Moreover, as in §7.2, set $t:=4bc-b-c=P\,\delta$.

Then the set of admissible triples $u=(\delta,b,c)^{\mathsf T}$ is described by the affine class $u_{0}+\Lambda$, where
$\Lambda\subset\mathbb{Z}^{3}$ is a sub-lattice of fixed index (independent of $P$), given by the linear congruences
\[
\boxed{\quad \delta^{\mathsf T} \equiv \delta\pmod{m_{3}},\qquad b^{\mathsf T}\equiv 0 \pmod{g},\qquad c^{\mathsf T}\equiv 0 \pmod{g}\quad}
\]
for any prefixed odd modulus $m_{3}$ satisfying $m_{3}\mid g$ and $\gcd(m_{3},P)=1$.

\paragraph{Proof of the linear conditions.}
From $b=g b'$ and $c=g c'$ it follows immediately that $b\equiv c\equiv 0\pmod{g}$, hence also $b\equiv c\equiv 0\pmod{m_{3}}$ when $m_{3}\mid g$.
Then
\[
t=4bc-b-c\equiv 0 \pmod{m_{3}},
\]
i.e., $P\,\delta \equiv 0\pmod{m_{3}}$. Since $\gcd(P,m_{3})=1$, we obtain $\delta \equiv 0\pmod{m_{3}}$.

On the other hand, from
\[
(4g b'-1)(4g c'-1)=4P\,\delta+1
\]
and the condition $m_{3}\mid g$ it follows that
\(
4g b_{1}-1\equiv -1 \pmod{m_{3}}
\)
and
\(
4g c_{1}-1\equiv -1 \pmod{m_{3}},
\)
therefore
\(
4P\,\delta+1\equiv (-1)\cdot(-1)\equiv 1 \pmod{m_{3}}
\),
i.e., $4P\,\delta\equiv 0\pmod{m_{3}}$, and since $\gcd(P,m_{3})=1$ we have $\delta\equiv 0\pmod{m_{3}}$.

Combining, we get $\delta^{\mathsf T}\equiv \delta \equiv 0\pmod{m_{3}}$, which is equivalent to the stated condition $\delta^{\mathsf T}\equiv \delta \pmod{m_{3}}$.
\paragraph{Vector form of the system.}
Let $u=(\delta,b,c)^{\mathsf T}$, then the system can be written as
\[
M\,u\equiv \rho \pmod m,\quad
M=\mathrm{diag}(1,1,1),\quad
\rho=\begin{pmatrix}\delta\\[2pt]0\\[2pt]0\end{pmatrix},\quad
m=\begin{pmatrix}m_{3}\\[2pt] g\\[2pt] g\end{pmatrix}.
\]
The set of all solutions is an affine lattice $u_{0}+\Lambda$, where
\[
u_{0}=\begin{pmatrix}\delta\\[2pt]0\\[2pt]0\end{pmatrix},\qquad
\Lambda=\{v\in\mathbb{Z}^{3}: v\equiv 0\ (\bmod\, (m_{3},g,g))\},
\]
its index is $[\mathbb{Z}^{3}:\Lambda]=m_{3}\,g^{2}$ and does not depend on $P$.
Note: when $m_{3}\mid g$ and $\delta\equiv 0\pmod{m_{3}}$ the first residue class is zero.

\paragraph{Numerical example \[(P=73, \alpha =1, d'=3) \].}
We have $g = \alpha  d' = 3$, $\delta = \alpha (d')^{2} = 9$.
Take $b_{1} = 1$, $c_{1} = 20$ (see~\S7.2), then $b = 3$, $c = 60$.
Choose modulus $m_{3} = 3$ such that $m_{3}\mid g$ and $\gcd(m_{3},P)=1$.
Check of the linear conditions: $b \equiv 0 \pmod{3}$, $c \equiv 0 \pmod{3}$, $\delta^{\mathsf T} \equiv \delta \equiv 9 \equiv 0 \pmod{3}$.
Additionally: $d' \mid (b_{1}+c_{1})$ and $\dfrac{b_{1}+c_{1}}{d'} = \dfrac{21}{3} = 7 \equiv 3 \pmod{4}$.
Thus $(\delta,b,c) = (9,3,60)$ lies in the affine class $u_{0} + \Lambda_{\mathrm{ED2}}$.


\subsection{Example: Type~I box for the parametrization from Section~\ref{sec:two-multiples} (ED2)}
\label{subsec:boxI_two_multiples}

In the ED2 subcase from §7.2 fix a square-free $\alpha$ and a number $d' \in \mathbb{N}$, setting
\[
  g := \alpha \,d', \qquad \delta := \alpha (d')^{2}.
\]
We consider triples $(A,B,C)$ with two multiples of $P$:
\[
  B = bP,\quad C = cP,\quad A \not\equiv 0 \pmod{P},
\]
where
\[
  b = g\,b',\quad c = g\,c',\quad \gcd(b',c') = 1.
\]

The ED2 linear conditions define an affine class.
\[
  u := (\delta,b,c)^{\mathsf T} \in u_0 + \Lambda_{\mathrm{ED2}}, \qquad
  u_0 := (\delta,0,0),
\]
where the sub-lattice
\[
  \Lambda_{\mathrm{ED2}} := \{\, v \in \mathbb{Z}^3 : v \equiv 0 \ (\bmod\ (m_3,\,g,\,g)) \,\},
\]
and the modulus $m_3$ is fixed, odd, divides $g$, and is co-prime to $P$.

Equivalently:
\[
  \delta^{\mathsf T} \equiv \delta \pmod{m_3}, \quad b \equiv 0 \pmod{g}, \quad c \equiv 0 \pmod{g}.
\]
The lattice index $[\mathbb{Z}^3 : \Lambda_{\mathrm{ED2}}] = m_3\,g^{2}$ does not depend on $P$.
\begin{lemma}
The index of the lattice \[ \Lambda_{\mathrm{ED2}} \] in $\mathbb{Z}^3$ is finite and equals
\[
[\mathbb{Z}^3 : \Lambda_{\mathrm{ED2}}] = m_3\,g^2
\]
where $m_3 \mid g$, $\gcd(m_3, P) = 1$. In particular, the set of points of the Type~I box has positive density.
\end{lemma}
\medskip
\noindent\textbf{Selection conditions in the Type~I box (linear + range):}
\begin{itemize}
  \item \textbf{Membership:} $u \in u_0 + \Lambda_{\mathrm{ED2}}$;
  \item \textbf{Order:} $b \le c$;
  \item \textbf{Range:} $1 \le \delta,b,c \le (\log P)^{\kappa}$ for fixed $\kappa>0$.
\end{itemize}

\noindent\textbf{Subsequent checks (outside the box):}
\begin{itemize}
  \item \textbf{GCD:} $\gcd(b,c) = g$ \; (equivalently $\gcd(b/g,\;c/g) = 1$);
  \item \textbf{Relation with $t$:} $t = 4bc - b - c = P\,d$;
  \item \textbf{Local conditions of §7.2:} $d'\mid (b'+c')$ and $\displaystyle \frac{b'+c'}{d'} \equiv 3 \pmod{4}$;
  \item \textbf{Integrality:} $A = \dfrac{bc}{\delta} \in \mathbb{N}$ and $A \not\equiv 0 \pmod{P}$.
\end{itemize}

\subsection{Box~II for parametrization of quadratic dependencies}
\label{subsec:boxII-quadratic}
In an earlier version, to describe surfaces above linear ones in the presence of quadratic dependencies, the construction of Box~II was introduced. In the current version these cases are handled differently, which allows leaving only Box~I in the main text.
\subsection{Method \mh{ED1} in the parametrization \mh{(\gamma, c, u, v)}}
\label{subsec:ED1-param}

The $ED1$ method is formulated in terms of the parameters $(\gamma, c, u, v)$ introduced in~\S\ref{sec:param-Bnot0}.
These parameters arise from the quadratic factorization equation
\[
  (\gamma A - c)(\gamma B - c) = c^{2},
\]
in which we set $u = \gamma A - c$ and $v = \gamma B - c$.

\paragraph{Note on the $ED1$ method.}
Unlike the $ED2$ case, the $ED1$ construction relies on the factorization identity
\[
  (\gamma A - c)(\gamma B - c) = c^{2},
\]
and after substituting $u = \gamma A - c$, $v = \gamma B - c$ leads to the nonlinear diophantine set
\[
  u v = c^{2}.
\]
This set is not an affine lattice in~$\mathbb{Z}^2$, so theorems on the density of points in affine lattices (see \S\ref{subsec:log-density}–\ref{subsec:log-convergence}) are not directly applicable.
In what follows, for~$ED1$ we use a combinatorial analysis over the divisors of~$c^{2}$.

An admissible quadruple $(\gamma, c, u, v)$ satisfies the conditions.
\[
  u v = c^{2}, \quad
  u \equiv v \equiv -c \pmod{\gamma}, \quad
  u \not\equiv -c \pmod{P}, \quad
  u \le v,
\]
where $\gamma$ and $c$ are related to~$P$ by $4c - 1 = \gamma P$, $\gcd(\gamma,c) = 1$.

\begin{lemma}[Counting admissible pairs for $ED1$]
\label{lem:ED1-div}
Let $m, n \ge 1$ be moduli and $a \bmod m$, $b \bmod n$ be residue classes. The number of pairs $(u,v) \in \mathbb{Z}_{>0}^2$
satisfying
\[
  u v = c^{2}, \quad u \equiv a \pmod{m}, \quad v \equiv b \pmod{n},
\]
equals the number of divisors $u \mid c^{2}$ lying in the intersection of two residue classes:
\[
  u \equiv a \pmod{m}, \quad u \equiv \overline{b} \, c^{2} \pmod{n},
\]
when the inverse element $\overline{b}$ exists (i.e. $(b, n) = 1$), and is zero otherwise.
In particular, if the system of congruences modulo $m$ and~$n$ is consistent modulo $M = \operatorname{lcm}(m,n)$, then
\[
  \#\{(u,v)\} \le \tau(c^{2}),
\]
and if there is at least one divisor $u \mid c^{2}$ in the given class modulo~$M$, the set is nonempty.
\end{lemma}

\paragraph{Summary for $ED1$.}
We do not use Theorems \S\ref{subsec:log-density}–\ref{subsec:log-convergence} for~$ED1$.
The existence and estimate of the number of admissible solutions are provided by Lemma~\ref{lem:ED1-div}
through a combinatorial analysis of the divisors of~$c^{2}$ and checking the consistency of linear congruences
for the divisor~$u$. Accordingly, all “density” conclusions of Section~\ref{sec:conv-completeness}
apply to the $ED2$ case, where the set of candidates is given by a system of linear congruences
(an affine lattice).
\begin{remark}
The density statements and conclusions of this section apply only to the set of $ED2$ solutions, which, unlike $ED1$, forms an affine finite-index lattice. For~$ED1$ there are no such lattice structures.
\end{remark}

\subsection{Method \mh{ED2}: case of two multiples of~\mh{P} \texorpdfstring{$(\delta,b,c)$}{(δ,b,c)}}
\label{subsec:ED2-two-mult}

Consider
\[
  \frac{4}{P} = \frac{1}{A} + \frac{1}{bP} + \frac{1}{cP}, \quad
  b, c \in \mathbb{N}, \quad
  A \not\equiv 0 \pmod{P}, \quad
  A \le bP \le cP.
\]
Multiplying by $A b c P$, we obtain the following.
\[
  A(4bc - b - c) = P b c.
\]
Set
\[
  t := 4bc - b - c = P \delta, \qquad \delta \in \mathbb{N},
\]
then
\[
  A = \frac{b c}{\delta}.
\]
The factorization
\[
  (4b - 1)(4c - 1) = 4P\delta + 1
\]
defines the parameters $r = 4b - 1$ and $s = 4c - 1$.

\begin{theorem}\label{thm:ED2-two-mult}
Let $P$ be prime, $\delta \in \mathbb{N}$ and
\[
  r s = 4P\delta + 1, \quad r \equiv s \equiv 3 \pmod{4}.
\]
Setting $b = \frac{r+1}{4}$, $c = \frac{s+1}{4}$, we obtain a solution
\[
  A = \frac{b c}{\delta}, \quad B = bP, \quad C = cP
\]
of the Erdős–Straus equation under the conditions
\[
  \delta \mid bc, \quad b \le c, \quad \frac{b c}{\delta} \le bP.
\]
These conditions are equivalent to the membership of $u$ in the affine lattice
\[
  u \equiv u_0(P) \pmod{\Lambda_1}
\]
with index~$M_1$ and additional local constraints
(primitivity, range $1 \le b, c, \delta \le (\log P)^A$, etc.).
\end{theorem}

\noindent
Denote by
\[
  \mathcal{C}_{\mathrm{ED2}}(P)
  = \left\{ (\delta,b,c) \in \mathbb{N}^3 \ \middle|\ 
      (4b - 1)(4c - 1) = 4P\delta + 1,\ 
      b \le c,\ 
      \delta \mid bc,\ 
      \tfrac{b c}{\delta} \le bP
  \right\}
\]
the set of all admissible parameter triples for fixed~$P$.
\subsection{Unconditionality for the \mh{ED2} algorithm}
\label{subsec:ED2-unconditional}
For any prime $P\equiv 1\pmod 4$, the geometric part of the $ED2$ method guarantees the existence of $(\delta,b,c)$ satisfying conditions (I)–(II) and the linear congruences, with $(\delta,b,c)$ lying in a Type~I parametric box. The transition to the full decomposition $\frac4P = \frac{1}{A}+\frac{1}{bP}+\frac{1}{cP}$ for fixed $P$ is completed \emph{without factorization} by normalizing the coordinates $(u,v)$ and the inverse point test (see \S D: Lemmas~\ref{lem:uv-square}, \ref{lem:back}): from $u=md'$, $u\equiv v\ (\mathrm{mod}\ 2)$ and $u^2-v^2=4A/\alpha$ we obtain $b'=(u+v)/2$, $c'=(u-v)/2$ and the required $b,c,A$.
\begin{remark}
The parametrization identity
\(
(4\alpha d' b'-1)(4\alpha d' c'-1)=4\alpha P{d'}^2+1
\)
holds automatically for the constructed $b',c'$ and can be used as a verification/algorithmic tool, but is not required for the proof of existence for fixed $P$.
\end{remark}

\begin{theorem}\label{thm:ED2-unconditional}
For any prime $P \equiv 1 \pmod{4}$ there exists a representation
\[
  \frac{4}{P} = \frac{1}{A} + \frac{1}{bP} + \frac{1}{cP},
\]
where $(\delta,b,c) \in \mathbb{N}^3$ are obtained from the constructive $ED2$ method, based on the parametrization of the set~$\mathcal{C}_{\mathrm{ED2}}(P)$ satisfying conditions~(I)–(III) below.
\end{theorem}

\emph{(I) Mathematical conditions:}
\begin{align*}
& b, c, \delta \in \mathbb{N}, \quad \gcd(b,c) = d\\
& (4b - 1)(4c - 1) = 4P\delta + 1, \\
& \delta \mid bc, \quad A = \frac{bc}{\delta} \le bP, \quad B = bP,\ C = cP, \\
& g_b = \gcd(b, g), \quad g_c = \gcd(c, g), \\
& b' = \frac{b}{g_b}, \quad c' = \frac{c}{g_c}, \quad \gcd(b',g) = \gcd(c',g) = 1, \\
& \alpha = \gcd\!\bigl(g,\, b'+c'\bigr),\quad d' = \frac{g}{\alpha}.
\end{align*}

Here $g$ and $d'$ are consistent with the parametrization of Theorem~\ref{thm:ED2-param}: $\delta = \alpha \,{d'}^2$, where $\alpha$ is square‑free.

\emph{(II) Algorithmic conditions (Type~I parametric box):}
\begin{itemize}
  \item $1 \le \delta,b,c \le (\log P)^{A_0}$, with fixed $A_0>0$.
  \item $b \le c$.
  \item $(\delta,b,c) \equiv u_0(P) \pmod{\Lambda_1}$, where $\Lambda_1 \subset \mathbb{Z}^3$ is a lattice of index $M_1$.
  \item $\gcd(b',c') = 1$.
  \item For the Type~I box in~$ED2$ there are no additional restrictions on parity or coprimality with~$P$.
\end{itemize}

\emph{(III) Unconditional guarantee of finding a solution.}

Consider the two‑dimensional lattice
\[
  L = \{\, (u,v)\in\mathbb{Z}^2 : u\,b' + v\,c' \equiv 0 \pmod{g} \,\}, 
\]
where $g\in\mathbb{N}$, $b',c'$ satisfy $\gcd(b',g)=\gcd(c',g)=1$, and set $\alpha := \gcd(g,\,b'+c')$, $d' := g/\alpha$.

\begin{lemma}[Kernel structure and diagonal period]\label{lem:alpha-step}
$L$ is a full‑rank lattice of index $g$, containing $\mathbf{v}_1=(c',-b')$ and $\mathbf{v}_2=(d',d')$. The vector $\mathbf{v}_2$ is a diagonal period of length $d'=g/\alpha$.
\end{lemma}

\begin{lemma}[Unique representative of a class]\label{lem:res-class}
For any $m\in\mathbb{N}$ and $r\in\mathbb{Z}$, in any half‑interval $[x_0,x_0+H)$ with $H\ge m$ there is exactly one integer $u$ with $u\equiv r\pmod m$.
\end{lemma}

\begin{lemma}[Diagonal coset]\label{lem:diag-coset}
If $\mathbf{w}=(d',d')\in L$ and $\mathbf{p}_0\in L$, then
\[
  \mathbf{p}_0+\mathbb{Z}\mathbf{w}=\{(u,v)\in L:\ u\equiv u_0,\ v\equiv v_0\pmod{d'}\}.
\]
\end{lemma}

\begin{proposition}[Hitting the box with a point of $L$]\label{prop:hit-box}
Let $R=[x_0,x_0+H)\times[y_0,y_0+W)\subset\mathbb{R}^2$ be an axis‑parallel rectangle (see the definition of the Type~I parametric box, \ref{def:boxI}). If $H\ge d'$ and $W\ge d'$, where $d'=d/\alpha$ from Lemma~\ref{lem:alpha-step}, then $L\cap R\neq\varnothing$.
\end{proposition}

\begin{proof}
Choose any point $\mathbf{p}_0=(u_0,v_0)\in L$ (for example, $\mathbf{v}_1=(c',-b')$ from Lemma~\ref{lem:alpha-step}). By Lemma~\ref{lem:res-class} there exist unique
\[
u^\ast\in[x_0,x_0+H)\ \text{with}\ u^\ast\equiv u_0\pmod{d'},\qquad
v^\ast\in[y_0,y_0+W)\ \text{with}\ v^\ast\equiv v_0\pmod{d'}.
\]
Set $m:=(u^\ast-u_0)/d'\in\mathbb{Z}$ and consider $\mathbf{w}=(d',d')$ from Lemma~\ref{lem:alpha-step}. Then the point
\[
\mathbf{p}:=\mathbf{p}_0+m\mathbf{w}=(u_0+md',\,v_0+md')
\]
belongs to $L$. By the choice of $m$ we have the first coordinate $u_0+md'=u^\ast\in[x_0,x_0+H)$. By Lemma~\ref{lem:diag-coset} the second coordinate of $\mathbf{p}$ lies in the same class modulo $d'$ as $v_0$, hence by the uniqueness of the class representative (Lemma~\ref{lem:res-class}) we get $v_0+md'=v^\ast\in[y_0,y_0+W)$. Therefore $\mathbf{p}\in L\cap R$.
\end{proof}

\begin{corollary}[Role of the parameter $\alpha$]\label{cor:alpha-role}
Increasing $\alpha$ decreases the diagonal step $d'=g/\alpha$ and thus relaxes the condition $H,W\ge d'$ in Proposition~\ref{prop:hit-box}. Equivalently, the density of the diagonal layers of $L$ in projection increases by a factor of $\alpha$. With the standard choice of Type~I box sizes (see~\ref{subsec:boxI}) the condition $H,W\ge d'$ is satisfied, and there exists an admissible point $L\cap R\neq\varnothing$.
\end{corollary}
\medskip\

\textbf{Connection with \S7.3.} Lemma~\ref{lem:alpha-step} and Proposition~\ref{prop:hit-box} provide a point $(b,c)$ satisfying the linear constraints of the affine class.

For this point to produce a solution to the ED2 problem, we use the normalization $(u,v)$ and the inverse test (Lemmas~\ref{lem:uv-square}, \ref{lem:back}): for $m=4A-P$ we take $u=md'$ and find $v\equiv u\ (\mathrm{mod}\ 2)$ with $u^2-v^2=4A/\alpha$.
\noindent
Combining Lemma~\ref{lem:alpha-step} and Proposition~\ref{prop:hit-box}, we obtain item~(III) of Theorem~\ref{thm:ED2-unconditional}:

\emph{reducing the lattice step via the parameter $\alpha$} gives an unconditional guarantee of the existence of an admissible triple $(\delta,b,c)\in\mathcal{C}_{\mathrm{ED2}}(P)$ in the given parametric box.

\begin{remark}
The conditions obtained in Theorem~\ref{thm:ED2-unconditional} have a natural algebraic interpretation in terms of the ED2 model. 
In particular, the parameters $m$, $M$, $A$ introduced in the geometric formulation correspond to the coefficients and constraints in the system~(ED2), 
where checking the existence of a solution reduces to analyzing congruences and inequalities. 
A detailed derivation, as well as extended criteria allowing refinement of the applicability bounds of the theorem, are given in Appendix~D. 
There it is also shown how the geometric construction of the window and strip is consistent with the algebraic description, 
and additional lemmas are provided for deeper consideration.
\end{remark}
\subsection{Key points addressing the factorization claim}\label{sec:no-factor}

\subsection*{What exactly is used in ED2}
- integer arithmetic and $\gcd$ operations;
- congruences modulo a prime $P$; computation of \(\Legendre{a}{P}\) and, if necessary, \(\Jacobi{a}{\gamma}\) — without factorization of the modulus;
- parity checks, conditions $u\equiv v\pmod 2$, $\gcd(u,v)=1$, equality $uv=c^2$.

\subsection*{What exactly is not used}
- factorization of $C$, $\gamma$ or intermediate numbers;
- root finding modulo composite moduli;
- factorization for the square test: the fact $uv=c^2$ is guaranteed by normalization.

We rely on:
- \textbf{Sum and discriminant} (\Cref{lem:uv-square})
- \textbf{Back‑test} (\Cref{lem:back})
- \textbf{Quadratic reparametrization} (\Cref{thm:quadratic})

\subsection*{Correctness via minimal lemmas}
We rely on:
- \textbf{Sum and discriminant} (\Cref{lem:uv-square}): $S=u$, $\Delta=v^2$;
- \textbf{Back‑test} (\Cref{lem:backtest}): conditions on $(u,v)$ for fixed $P$;
- \textbf{Quadratic reparametrization} (\Cref{thm:quadratic}).

\begin{proposition}[ED2 does not rely on factorization]\label{prop:no-factor}
The ED2 algorithm, forming and checking pairs $(u,v)$ for fixed prime $P$, uses only:
(i) integer operations and $\gcd$; (ii) checks modulo $P$ and the symbols \(\Legendre{\cdot}{P}\)/\(\Jacobi{\cdot}{\gamma}\);
(iii) linear formulas for recovering $A,B$ and setting $C=cP$ from $uv=c^2$.
At no step of ED2 is factorization required.
\end{proposition}
\begin{proofsketch}
\Cref{lem:uv-square} gives the link to $(u,v)$; \Cref{lem:backtest} — sufficiency of local conditions; $C=cP$ follows from $uv=c^2$. All checks reduce to arithmetic, $\gcd$, and Legendre/Jacobi symbols.
\end{proofsketch}

\begin{remark}[Optional accelerations]
Sieving by small primes is permissible as an optimization, but is not part of the correctness of ED2 and is not used in \Cref{lem:uv-square,lem:backtest,thm:quadratic,prop:no-factor}.
\end{remark}
\subsection{Enumeration algorithm for the~\mh{ED2} method}
\label{subsec:ED2-enum}

As follows from Theorem~\ref{thm:ED2-unconditional}, an admissible triple $(\delta,b,c)$ exists for any prime $P \equiv 1 \pmod{4}$.


\[ \mathcal{D} := \{\, \delta \le X : \delta \equiv 3 \ (\mathrm{mod}\ 4) \,\}, \qquad M(\delta):=P+\delta.\]

\medskip
\textbf{Analytic part (for BV).}
For each $\delta \in \mathcal{D}$ set
\[
S(\delta) := \#\{\, a \le X^2 : a \ \text{is prime},\ a \equiv -1 \pmod{\delta}\,\}.
\]
On average over $\delta \le X$ and for $\delta \le X/(\log X)^B$ (arbitrarily fixed $B>0$) there holds a Bombieri–Vinogradov type estimate:
\[
\sum_{\substack{\delta \le X \\ \delta \equiv 3 \ (\mathrm{mod}\ 4)}}
\Bigl|\, S(\delta) - \frac{X^2}{\varphi(\delta)\,\log X} \Bigr|
\ \ll_{B}\ \frac{X^2}{(\log X)^{A}},
\]
for any fixed $A>0$.
\begin{lemma}[Parameterization via $(t,k)$]
\label{lem:param-tk-inline}
There exist natural numbers $\delta,a$ with $a \mid (P+\delta)$ and $a \equiv -1 \pmod{\delta}$ if and only if
there exist $t,k \in \mathbb{N}$ such that, with $D:=tk-1$, one has $D \mid (P+t)$ and $D \mid (kP+1)$, and
\[
\delta = \frac{P+t}{D}, \qquad a = \frac{kP+1}{D}.
\]
In this case $P+\delta = t\,a$ and $a = k\delta-1$.
\end{lemma}
\medskip
\textbf{Algorithmic part (with divisibility).}
Define
\[
U(\delta) := \#\{\, a \le X^2 : a \ \text{is prime},\ a \mid M(\delta),\ a \equiv -1 \pmod{\delta}\,\},
\qquad
I(\delta) := \mathbf{1}_{\{U(\delta)\ge 1\}}.
\]
To ensure $I(\delta)=1$ we use the parameterization via $(t,k)$ from Lemma~\ref{lem:param-tk-inline}:
if $D=tk-1$ divides both $(P+t)$ and $(kP+1)$, then
\[
\delta = \frac{P+t}{D}, \qquad a = \frac{kP+1}{D},
\]
and then automatically $a \mid M(\delta)$ and $a \equiv -1 \pmod{\delta}$.

\textbf{Constructive procedure for finding such a solution}.

For a fixed prime~$P$:
\begin{enumerate}
  \item \textbf{Loop over $\delta$:} choose $\delta$ in the search range.
  \item \textbf{Factorization:} set \(N_\delta := 4P\delta + 1\) and factor \(N_\delta = r \cdot s\) with \(r \equiv s \equiv 3 \pmod{4}\).
  \item \textbf{Pass to $(b,c)$:} $b = (r+1)/4$, $c = (s+1)/4$.
  \item \textbf{Filters:} $\delta \mid bc$, $b \le c$, $bc/\delta \le bP$.
  \item \textbf{Construction:} $A = bc/\delta$, $B = bP$, $C = cP$.
  \item \textbf{Normalization:} remove duplicates arising from swapping $b \leftrightarrow c$ in symmetric cases.
\end{enumerate}

\noindent
Enumerating the elements of $\mathcal{C}_{\mathrm{ED2}}(P)$ with filters~(I)–(II) from Theorems~\ref{th:log_density}–\ref{th:log_convergence}
guarantees finding at least one admissible solution by the~$ED2$ method in time polynomial in~$\log P$.
In particular, for $T = (\log P)^{A_0}$, the set of tested triples has a positive density
in the corresponding class $u_0(P) \pmod{\Lambda_1}$, which ensures the finiteness of the search and no omissions.

\emph{Caution.}
Statements about the running time polynomial in $\log P$ are valid in the presence of a factorization oracle
(or conditionally - under standard heuristics for the average factorization time of \(N_\delta\)).


\subsubsection{Correspondence between~\mh{ED2} and~\mh{ED1} in new parameters}
\label{subsubsec:ED2-ED1-map}

The~$ED2$ and~$ED1$ methods use different parameterizations but produce solutions to the same Erdős–Straus equation for $P \equiv 1 \pmod{4}$.
The parameters $(\delta, b, c)$ in~$ED2$ and $(\gamma, c, u, v)$ in~$ED1$ are related constructively.

\begin{theorem}
\label{thm:ED2-to-ED1}
Let $(\delta, b, c) \in \mathcal{C}_{\mathrm{ED2}}(P)$ be an admissible~$ED2$ triple, and let
\[
  A = \frac{b c}{\delta}, \quad B = bP, \quad C = cP.
\]
Set
\[
  \gamma := \frac{4c - 1}{P}, \quad u := \gamma A - c, \quad v := \gamma B - c.
\]
Then $(\gamma, c, u, v)$ is an admissible~$ED1$ quadruple satisfying:
\begin{align*}
& (\gamma A - c)(\gamma B - c) = c^2, \\
& u \mid c^2, \quad u \equiv v \equiv -c \pmod{\gamma}, \quad u \le v, \quad u \not\equiv -c \pmod{P}, \\
& \gcd(\gamma, c) = 1, \quad \gamma \equiv 3 \pmod{4}.
\end{align*}
\end{theorem}

\begin{proof}[Idea]
From the factorization $(4b - 1)(4c - 1) = 4P\delta + 1$ it follows that $t := 4bc - b - c = P\delta$.
Substituting into the~$ED2$ equation we get $A = bc/\delta$, $B = bP$, $C = cP$.
The transition to~$ED1$ is made via
\[
  \gamma = \frac{4c - 1}{P}, \quad u = \gamma A - c, \quad v = \gamma B - c.
\]
Then
\[
  (\gamma A - c)(\gamma B - c) = \gamma^2 AB - \gamma c(A + B) + c^2 = c^2,
\]
and the residues $u \equiv v \equiv -c \pmod{\gamma}$ are obvious. The condition $u \mid c^2$ follows from the structure
$u = \gamma A - c = \tfrac{\gamma b c - c\delta}{\delta}$ and the equality $t = P\delta$;
$\gcd(\gamma, c) = 1$ comes from $4c - 1 = \gamma P$.
Since $4c - 1 \equiv 3 \pmod{4}$, we have $\gamma \equiv 3 \pmod{4}$.
\end{proof}

\medskip
\noindent
Similarly, for any admissible quadruple $(\gamma, c, u, v)$ of the~$ED1$ method
one can recover the triple $(\delta, b, c)$ of the~$ED2$ method via:
\[
  A := \frac{u + c}{\gamma}, \quad
  B := \frac{v + c}{\gamma}, \quad
  b := \frac{B}{P} = \frac{v + c}{\gamma P}, \quad
  \delta := \frac{b c}{A}.
\]
Verification of the~$ED2$ conditions is carried out according to the definition of the set~$\mathcal{C}_{\mathrm{ED2}}(P)$.

\subsubsection{Convolution: transition from~\mh{ED2} to~\mh{ED1}}
\label{subsubsec:ED2-to-ED1}

\paragraph{Introduction and motivation.}
In what follows, we will need to establish a connection between solutions of types~$ED2$ and~$ED1$.
Recall that the set~$\mathcal{C}_{\mathrm{ED1}}(P)$ of admissible quadruples $(y, c, u, v)$ for the~$ED1$ method
is given in Definition~\ref{def:ED1_set_new} and includes divisibility conditions, congruences, and the identity $uv = c^2$.

\medskip
\noindent\textbf{Why “convolution” is needed.}
The procedure described below allows one to obtain from a correct~$ED2$ solution for a prime~$P$
a correct~$ED1$ solution (for the parameter $P' = \tfrac{4c - 1}{y}$) belonging to an important subclass
in which \emph{exactly one denominator} is divisible by~$P$.
Direct construction of this subclass within~$ED1$ is associated with serious technical difficulties,
therefore “convolution” serves as a constructive method for obtaining such solutions and as a tool for transferring structural properties between the methods.
This is not a one-to-one correspondence between all~$ED2$ and~$ED1$ solutions for a fixed~$P$,
but an auxiliary construction for classification purposes.

\paragraph{Admissible parameters of the~$ED1$ method.}
\begin{definition}[Admissible parameters of $ED1$]
\label{def:ED1_set_new}
Let $P$ be a prime, $P \equiv 1 \pmod{4}$. A quadruple $(y, c, u, v)$ belongs to $\mathcal{C}_{\mathrm{ED1}}(P)$ if:
\begin{enumerate}
  \item $y, c, u, v \in \mathbb{N}$, \quad $u \le v$;
  \item $y \mid (4c - 1)$, \quad $y \equiv 3 \pmod{4}$;
  \item $\gcd(y, c) = 1$;
  \item $u \mid c^2$;
  \item $v = \dfrac{c^2}{u}$.
\end{enumerate}
\end{definition}

\paragraph{Convolution algorithm $\boldsymbol{ED2 \to ED1}$.}
Let $(A, B, C, \delta, b, c) \in \mathcal{C}_{\mathrm{ED2}}(P)$, where $P$ is a prime, $P \equiv 1 \pmod{4}$. Then:
\begin{enumerate}
  \item $c \coloneqq C/P$, \quad $s \coloneqq 4c - 1$;
  \item choose the minimal $y \mid s$ with $y \equiv 3 \pmod{4}$;
\item (Canonical choice) Set \(y:=\gamma=(4c-1)/P\). Then \(P'':= (4c-1)/y=P\), and the subsequent steps
yield a correct~$ED1$ quadruple for \emph{the same} prime \(P\).

  \item $u \coloneqq yA - c$, \quad $v \coloneqq c^2 / u$;
  \item construct
    \[
      A' = A, \quad
      B' = \frac{u + v}{P'}, \quad
      C' = \frac{uv + 1}{P'}.
    \]
\end{enumerate}
These formulas are consistent with the general transition described in~\S9.9.1.
\begin{remark}
If one takes an arbitrary divisor \(y \mid (4c-1)\), then \(P'' := (4c-1)/y\) is in general unrelated to the original \(P\)
(the equality \(P=4P''+1\) does not hold). The transition to~$ED1$ is then correct only under the additional requirement “\(P''\) is prime”,
and the result refers to modulus \(P''\), not to the original \(P\).
\end{remark}

\begin{theorem}[Correctness of convolution]
Let \(P \equiv 1 \ (\mathrm{mod}\ 4)\) be a prime and \((A,B,C,\delta,b,c) \in \mathcal{C}_{\mathrm{ED2}}(P)\).
Set \(\gamma := (4c-1)/P\), \(y:=\gamma\), \(u:=yA-c\), \(v:=c^2/u\).
Then \((y,c,u,v) \in \mathcal{C}_{\mathrm{ED1}}(P)\) and
\[
(yA-c)(yB-c)=c^2,\quad u\equiv v \equiv -c \pmod y,\quad u \mid c^2,\quad u \le v,\quad \gcd(y,c)=1,\quad y \equiv 3 \ (\mathrm{mod}\ 4).
\]
\end{theorem}
\begin{proof}
By setting \(y:=\gamma=(4c-1)/P\) (canonical choice), we obtain the statement as a special case of Theorem~\ref{thm:ED2-to-ED1}.
The conditions \(u\equiv v\equiv -c \pmod y\), \(u\mid c^2\), \(\gcd(y,c)=1\), \(y\equiv 3\pmod 4\) follow from the definition of \(\gamma\) and the transition formulas.
\end{proof}

\begin{proposition}[Characterization of the image]
\label{prop:image}
Let $(\gamma, c, u, v)$ be an admissible~$ED1$ solution,
$A = (u + c)/\gamma$, $B = (v + c)/\gamma$.
Then $(\gamma, c, u, v)$ lies in the image of the convolution if and only if
$(A, B)$ satisfies the linear–modular admissibility constraints of~$ED2$.
\end{proposition}

\begin{corollary}[Conditional completeness]
\label{cor:conditional-surj}
If the~$ED2$ enumeration covers all affine classes modulo
$M = \mathrm{lcm}(P, \gamma)$ induced by admissible~$ED1$ solutions,
then the convolution (see \ref{def:ED1_set_new}) is surjective onto the set of~$ED1$ solutions.
\end{corollary}

\paragraph{Implementation note.}
In practice, the “full coverage” condition is replaced by a multi-coset enumeration of~$ED2$ classes
sufficient to cover all compatible classes.
When the condition of Corollary~\ref{cor:conditional-surj} is met, the convolution becomes complete.

\begin{corollary}
\label{cor:ED2_to_ED1_power_newparams}
If $\mathcal{C}_{\mathrm{ED2}}(P) \neq \varnothing$ and there exists
$(A, B, C, \delta, b, c) \in \mathcal{C}_{\mathrm{ED2}}(P)$
for which $u = yA - c$ and $v = c^2/u$ admit at least two distinct choices of
$y \mid (4c - 1)$, $y \equiv 3 \pmod{4}$,
yielding different $P'$, then
\[
  \bigl|\mathcal{C}_{\mathrm{ED1}}^{(\mathrm{conv})}(P)\bigr|
    > \bigl|\mathcal{C}_{\mathrm{ED2}}(P)\bigr|.
\]
\end{corollary}

\begin{remark}
The example $P = 2521$, in which three~$ED2$ solutions generate six~$ED1$ solutions,
illustrates the accuracy of the mapping and the potential completeness under full coverage of the coset,
but does not prove surjectivity in the general case.
\end{remark}
\subsection{Anticonvolution for \mh{ED1 - ED2} and the Canon Existence Hypothesis}
\label{subsec:ED1-ED2-antisvertka}
\begin{figure}[h!]
\centering
\begin{tikzpicture}[
  node distance=12mm,
  every node/.style={font=\small},
  ed2/.style={rectangle,draw,rounded corners,fill=blue!10,inner sep=4pt},
  ed1/.style={rectangle,draw,rounded corners,fill=green!10,inner sep=4pt},
  arrow/.style={-{Latex[length=2mm]},thick}
]

\node[ed2] (ed2) {ED2 solution};

\node[ed1a,ed1,above right=of ed2,xshift=14mm] (ed1a) {ED1: $(y,c,u_1,v_1)$};

\node[ed1b,ed1,below right=of ed2,xshift=14mm] (ed1b) {ED1: $(y,c,u_2,v_2)$};

\node[ed1c,ed1,below=of ed1b] (ed1c) {ED1: $(y,c,u_3,v_3)$};

\draw[arrow] (ed2) -- (ed1a) node[midway,above left]{\small anticonvolution};
\draw[arrow] (ed2) -- (ed1b);
\draw[arrow] (ed2) -- (ed1c);

\draw[arrow,dashed] (ed1a.west) .. controls +(left:15mm) and +(up:10mm) .. (ed2.north west)
  node[pos=0.3,left]{\scriptsize $B \mid P$};

\node[red!80!black] at ($(ed1b.east)+(6mm,0)$) {\Large$\times$};
\node[red!80!black] at ($(ed1c.east)+(6mm,0)$) {\Large$\times$};

\end{tikzpicture}
\caption{Branching scheme in anticonvolution ED2$\rightarrow$ED1 and invertibility conditions}
\label{fig:2}
\end{figure}
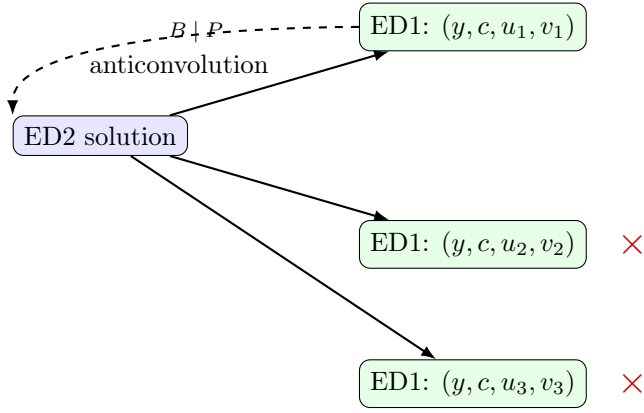

\noindent
\textbf{On an intuitive level:}
From any type~2 solution one can reconstruct at least one branch of type~1 solutions [fig.\ref{fig:2}], and then an entire “fan” of solutions by varying the parameter $u$ within the same residue class.
Branching occurs because, for fixed $(y,c)$ in ED1, several values of $u$ are admissible, each giving its own pair $(A,B)$ with the same long denominator $C = cP$.
In the reverse direction (ED1~$\rightarrow$~ED2) the transition is possible only for those branches where $B$ is divisible by $P$.

\medskip

\subsubsection*{Introduction: canonical subclass and completeness of solutions}
\paragraph{Hypothesis.}
In transitions between the $ED1$ and~$ED2$ representations, two main problems arise:
\emph{multiplicity} of admissible parameters for the same solution
and \emph{incompleteness} in the reverse transition.
To eliminate them, a \emph{canonical subclass} is introduced — a subset of parameters
for which the \emph{convolution} and \emph{anticonvolution} procedures are mutually one-to-one.

\medskip
\noindent\textbf{Motivation.}
Within the canon, each $ED1$ configuration has exactly one image in~$ED2$ and vice versa,
which ensures \emph{completeness} of the set of solutions and removes ambiguities
associated with permutations, scaling, and parameter choices.

\paragraph{Canonization parameters.}
In the $ED2$ representation, integers $m, o \ge 1$ are given,
serving as canonization moduli, and integers $d, b, c \in \mathbb{Z}$.
The parameter $d$ is defined as
\[
d \equiv y^{-1} \pmod{m},
\]
where $y$ is the $ED1$ quadruple component and $(y,m) = 1$.
The element $b$ is related to $u$ by the condition
\[
b \equiv u \pmod{o},
\]
and $c$ coincides in both the $ED1$ and $ED2$ structures.

\paragraph{Definition of the canonical subclass.}\label{par:}
An $ED1$ quadruple $\langle y, c, u, v \rangle$\label{eq:ED1_canon} and an $ED2$ triple $\langle d, b, c \rangle$\label{eq:ED2_canon}
belong to the \emph{canonical subclass} if:
\[
\gcd(y, c) = 1, \quad \gcd(u,v) = 1, \quad u < v, \quad uv = c^2,
\]
\[
u \equiv -c \pmod{y}, \quad b \equiv u \pmod{o}, \quad
0 < c < \min(m, o), \quad d \equiv y^{-1} \pmod{m}.
\]
\begin{lemma}[Anticonvolution formula]  
\label{lem:anti}  
Let $m,o$ be natural numbers in the canon with $(m,o)=1$, and let $(y,mo)=1$.  
Let $d$ denote the inverse of $y$ modulo $mo$, i.e., $d \equiv y^{-1} \pmod{mo}$.  
Then for $A = \dfrac{u+c}{y}$ the congruence holds  
\begin{equation}  
\label{eq:anti-master}  
A \equiv d\,(u+c) \pmod{mo}.  
\end{equation}  
\end{lemma}  

\begin{proof}  
By definition of $A$ we have the identity $yA = u + c$ in the integers. Passing to residue classes modulo $mo$ and multiplying by $d \equiv y^{-1} \pmod{mo}$, we obtain   

$A \equiv d\,(u+c) \pmod{mo}$

This is exactly \eqref{eq:anti-master}.  
\end{proof}  

\begin{remark}[On “reading Bézout’s lemma backwards”]
\label{rem:bezout-backwards}
Sometimes there is a temptation to “read Bézout’s lemma backwards” and immediately obtain congruences for the variable \(A\).
The correct interpretation here is only the equivalence
\[
\gcd(y,M)=1 \;\Longleftrightarrow\; \exists\,p,q\in\mathbb{Z}:\; p\,y+q\,M=1,
\]
that is, the existence of the inverse element \(y^{-1}\) modulo \(M\). However, this equivalence by itself does \emph{not} impose a congruence on \(A\) without an additional link between \(A\) and \(y\).

\medskip
\noindent\textbf{Where the naive reasoning goes wrong.}
From \(d \equiv y^{-1} \pmod m\) it follows only that \(d\,y \equiv 1 \pmod m\).
The statement \(A \equiv d\,y \pmod m\) is equivalent to \(A \equiv 1 \pmod m\),
which in general is false: there is no information here linking \(A\) with \(y\) modulo \(m\).

\medskip
\noindent\textbf{Correct derivation of the “anticonvolution formula”.}
The key step is to use the structural identity from ED1:
\[
yA = u + c \quad\text{in } \mathbb{Z}
\quad\bigl(\text{equivalently } u \equiv -c \pmod y \text{ and } A=\tfrac{u+c}{y}\bigr).
\]
Let \(M\) be the gluing modulus; in the canon it is convenient to take \(M=m\,o\) with \(\gcd(m,o)=1\) and \(\gcd(y,M)=1\).
Then, from the existence of the inverse \(d \equiv y^{-1} \pmod M\), multiplying both sides of the congruence \(yA\equiv u+c\pmod M\) yields
\[
A \equiv d\,(u+c) \pmod M,
\]
which gives the correct “anticonvolution formula”.

\medskip
\noindent\textbf{Conclusion.}
Bézout’s lemma ensures only the invertibility of \(y\) modulo the chosen modulus; the congruence for \(A\)
arises \emph{only} after invoking the link \(yA=u+c\). Therefore, “reading Bézout’s lemma backwards”
without this link is incorrect.
\end{remark}
\paragraph{Comments.}
\begin{itemize}
  \item The conditions on $\langle y, c, u, v \rangle$ coincide with those stated above in this section for ~$ED1$ (\ref{eq:ED1_canon}).
  \item The conditions on $\langle d, b, c \rangle$ coincide with those stated above in this section for ~$ED2$ (\ref{eq:ED2_canon}).
  \item The formula applies only within the canon; outside it, mutual one-to-one correspondence is not guaranteed.
\end{itemize}

\begin{theorem}[Non-emptiness of $ED2$ from non-emptiness of $ED1$]
\label{th:ed1-to-ed2-nonempty}
If $P \equiv 1 \pmod{4}$ and $\mathcal{C}_{\mathrm{ED1}}(P) \neq \varnothing$, then anticonvolution
produces a non-empty set $\mathcal{C}_{\mathrm{ED2}}^{(\mathrm{anti})}(P)$.
\end{theorem}

\begin{proof}
For any $\langle y, c, u, v \rangle \in \mathcal{C}_{\mathrm{ED1}}(P)$, formula~\eqref{eq:anti-master}
and the $ED1/ED2$ canonization yield a unique $\langle d, b, c \rangle \in \mathcal{C}_{\mathrm{ED2}}(P)$.
\end{proof}

\begin{corollary}[Power compression]
\label{cor:power-compression}
It is possible that $\bigl|\mathcal{C}_{\mathrm{ED2}}^{(\mathrm{anti})}(P)\bigr| < |S|$
for $S \subset \mathcal{C}_{\mathrm{ED1}}(P)$, since different $ED1$ quadruples
may correspond to the same $ED2$ triple.
\end{corollary}

\begin{remark}
For $P = 2521$, several distinct $ED1$ quadruples via anticonvolution yield the same $ED2$ triple,
illustrating the power compression effect.
\end{remark}

\section{Experimental data [tab.\ref{tab:3}]}

\begin{table}[!ht]
\centering
\caption{Examples of the algorithm in the forward and reverse directions}
\begin{tabular}{|c|c|c|c|c|c|c|c|c|c|c|c|}
\hline
\multicolumn{6}{|c|}{ED1 $\rightarrow$ ED2} & \multicolumn{6}{c|}{ED2 $\rightarrow$ ED1} \\
\hline
$P$ & $A$ & $B$ & $C$ & $y$ & $c$ & $P'$ & $A'$ & $B'$ & $C'$ & Invariants & Success \\
\hline
97 & 34 & 85 & 16490 & 7 & 170 & 129 & 34 & 645 & 21930 & $A,c$ & \checkmark \\
97 & 26 & 364 & 35308 & 15 & 364 & 103 & 26 & 2884 & 37492 & $A,c$ & \checkmark \\
\hline
\end{tabular}
\label{tab:3}
\end{table}
\section{Conclusion}
The combined use of the geometric construction (\emph{hitting the diagonal layer into the window}; see \ref{prop:hit-box}) and the algebraic model ED2 (Appendix~B), as well as its geometric refinement in Appendix~D, yields the following.

- The general case: namely, when $P$ runs over infinite sets, is certainly not proven constructively in the present work; the sections of Appendix~D that use Dirichlet’s theorem and finite coverings are conditional in nature.

- Special case: for any fixed odd prime $P$, the existence of a solution is strictly proven. The algebraic and geometric conditions are consistent: By Lemma \ref{lem:uv-square} we have \(S=u\), \(\Delta=v^2\), and Proposition \ref{prop:hit-box} ensures that the diagonal layer hits the target window (see also Theorem~9.21).

\section*{Acknowledgements}
The author expresses sincere gratitude to various AI models for assistance in refining explanations, proofreading, and creating illustrative Python examples.

\section*{Funding statement}
No funding was received.

\section*{Conflict of interest}
We have no conflicts of interest to disclose.

\section*{ Data availability}
The data that support the findings of this study are available from the corresponding author, upon reasonable request.
\nocite{*}
\bibliographystyle{plain}
\bibliography{esc_bib}

\appendix

\section{Notation and analytical tools}
\label{sec:app-tools}

\begin{theorem}[Bombieri–Vinogradov, large sieve]
\label{th:bombieri-vinogradov}
For any $A > 0$ there exists $y_0 = y_0(A)$ such that
\[
  \sum_{q \le Q} \max_{\substack{(a,q)=1}}
    \left| \pi(y; q, a) - \frac{y}{\varphi(q)} \right|
    \ll \frac{y}{(\log y)^A}
\]
uniformly for $Q \le y^{1/2} / (\log y)^{C(A)}$ and $y \ge y_0$.
\end{theorem}

\begin{theorem}[Greaves, larger sieve]
\label{th:greaves-larger-sieve}
If $\mathcal{A} \subset \{1,\ldots,N\}$ excludes too many residue classes for many primes $p \le B$, then
\[
  |\mathcal{A}| \le N \exp\!\left( -c \,\frac{\log B}{\log\log B} \right),
\]
where $c > 0$ depends only on the proportion of excluded classes.
\end{theorem}

\begin{theorem}[Chebotarev]
\label{th:chebotarev}
Let $L/K$ be a finite normal extension of number fields with Galois group $G$, and let $C \subset G$ be a conjugacy class.
Denote by $\pi_C(x)$ the number of prime ideals $\mathfrak{p}$ in~$K$ with norm $N\mathfrak{p} \le x$ for which $\mathrm{Frob}_{\mathfrak{p}} \in C$.
Then
\[
  \pi_C(x) = \frac{|C|}{|G|} \,\mathrm{Li}(x)
    + O\!\left( x \,\exp\!\bigl(-c\,\sqrt{\log x}\,\bigr) \right),
\]
where $c > 0$ and the constant in the $O$‑term depend only on the extension $L/K$.
\end{theorem}

\begin{remark}
\label{rem:not-using-sieves-direct}
We do not use Theorems~\ref{th:bombieri-vinogradov}–\ref{th:greaves-larger-sieve} to extract divisors of logarithmic order from fixed numbers; their purpose is to relate the method to known heuristic arguments.
\end{remark}

\bigskip
\section{Local application}
\label{sec:app-local}
Let $X \ge 2$ and
\[
\mathcal{D}_X := \{\, \delta \le X : \delta \equiv 1,3 \ (\mathrm{mod}\ 4) \,\}.
\]
Define the \emph{analytic} counter
\begin{equation}\label{eq:td_b}
T(\delta;X) = \frac{X^2}{\varphi(\delta)\,\log X} + O\!\left(\frac{X^2}{(\log X)^{1+\eta}}\right).
\end{equation}
Separately, for algorithmic purposes, introduce the \emph{divisibility} counter
\[
U(\delta;X) := \#\{\, a \le X^2 : a \ \text{prime},\ a \mid (P+\delta),\ a \equiv -1 \pmod{\delta} \,\},
\qquad
I(\delta;X) := \mathbf{1}_{\{U(\delta;X)\ge 1\}}.
\]
In this appendix, we work only with $T(\delta;X)$; the divisibility condition $a \mid (P+\delta)$ will be ensured
constructively in,9.10 through the parameterization $(t,k)$.


From Theorem A.1 (Bombieri–Vinogradov form with averaging over moduli) it follows that for any fixed $A,B>0$
and all sufficiently large $X$ we have the following.
\[
\sum_{\substack{\delta \le X/(\log X)^B \\ \delta \equiv 1,3 \ (\mathrm{mod}\ 4)}}
\Bigl|\, T(\delta;X) - \frac{X^2}{\varphi(\delta)\,\log X} \Bigr|
\ \ll_{A,B}\ \frac{X^2}{(\log X)^{A}}.
\]

In particular, for any fixed $\eta\in(0,1)$  

\[
\#\Bigl\{\delta \le \tfrac{X}{(\log X)^B} : \delta \equiv 3 \pmod{4},\   
\bigl|T(\delta;X)-\tfrac{X^2}{\varphi(\delta)\log X}\bigr|>\tfrac{X^2}{(\log X)^{1+\eta}}\Bigr\}  
\ \ll_{A,B,\eta}\ \frac{X}{(\log X)^{A-\eta}}.  
\]  

Hence for most $\delta \le X/(\log X)^B$ we have

\[
T(\delta;X) = \frac{X^2}{\varphi(\delta)\,\log X} + O\!\left(\frac{X^2}{(\log X)^{1+\eta}}\right),
\]

where $\eta>0$ is a fixed constant, and the constants in the $O$‑terms are absolute.

Choosing $A-\eta > B$ (for example, $A=B+2$ and $\eta=\tfrac12$), we obtain that the exceptional set has
size $o\!\left(X/(\log X)^B\right)$, i.e., for a proportion $1-o(1)$ of moduli $\delta \le X/(\log X)^B$ relation [\ref{eq:td_b}] holds.

This estimate is the main asymptotic within the large sieve bound and is refined in Theorem B.1 below.

\begin{remark}[Do not mix with divisibility]
\label{rem:B-no-div}
The counter $T(\delta;X)$ does \emph{not} include the condition $a \mid (P+\delta)$ and is used only to estimate the density
of primes in the class $-1 \pmod{\delta}$. The divisibility condition will be ensured in §\,9.10 by the construction via pairs $(t,k)$,
which yields specific $(\delta,a)$ with $a \mid (P+\delta)$ and $a \equiv -1 \pmod{\delta}$, i.e., $U(\delta;X)\ge 1$.
\end{remark}

\section{Existence of \mh{ED2}-type solutions}
\label{app:existence-ED2}
\label{appC:bridge-ED2}

Let $P$ be an odd prime, $X:=P^{1/2}$, and
\[
\mathcal{D}_X := \{\, \delta \le X : \delta \equiv 3 \ (\mathrm{mod}\ 4) \,\},
\qquad
M(\delta):=P+\delta.
\]
We use two counters (analogous to Appendix~B):
\[
S(\delta;X) := \#\{\, a \le X^2 : a \ \text{prime},\ a \equiv -1 \pmod{\delta}\,\},
\]
\[
U(\delta;X) := \#\{\, a \le X^2 : a \ \text{prime},\ a \equiv -1 \pmod{\delta},\ a \mid M(\delta)\,\}.
\]
Bombieri–Vinogradov type estimates apply to $S(\delta;X)$ (see Appendix~B) and \emph{do not} apply to $U(\delta;X)$;
divisibility in $U$ is ensured constructively via §\,\ref{subsec:ED2-enum}.

\medskip
\noindent
Connection with §\,\ref{subsec:ED2-enum}.
Lemma~\ref{lem:param-tk-inline} gives a parameterization $(t,k)\mapsto(\delta,a)$ with
$a \mid M(\delta)$ and $a \equiv -1 \pmod{\delta}$; moreover, $4 \mid M(\delta)$, so $4\alpha\mid M(\delta)$.
The construction of the “raw” pair $(b,c)$ and its subsequent normalization for ED2 is carried out in §\,\ref{subsec:ED2-enum};
here we record their properties.

\begin{lemma}[Normalization of the pair]
\label{lem:C-normalization}
Let $\delta\in\mathcal{D}_X$ and suppose $a\le X^2$ and integers $b,c$ are found as described in §\,\ref{subsec:ED2-enum}
(in particular, $a\mid \gcd(b,c)$ and the ED2 local congruences hold for the “raw” pair $(b,c)$).
Set $d:=\gcd(b,c)$, $d':=d/a$, $b':=b/d'$, $c':=c/d'$.
Then
\[
\gcd(b',c')=1,
\quad
\text{the ED2 local conditions carry over from $(b,c)$ to $(b',c')$},
\]
and the crude bounds
\[
8\le \delta \le X,\qquad a\le X^2,\qquad b',c'\ \ll\ P^{3/2}\,(\log P)^{O(1)}
\]
hold unconditionally.
\end{lemma}

\begin{proof}[Proof sketch]
By the construction in §\,\ref{subsec:ED2-enum} we have $a\mid d$, hence $d'\in\mathbb{N}$ and $(b',c')$ are integers.
Any congruences and divisibility conditions tied to $b,c$ and the modulus $\delta$ are preserved upon division by $d'$.
The bounds follow from the ranges built into §\,\ref{subsec:ED2-enum}.
\end{proof}

\begin{lemma}[Finiteness of the enumeration]
\label{lem:C-termination}
Fix an explicit upper bound on the moduli $\Delta_{\max}$ (for example, $\Delta_{\max}=\lfloor X\rfloor$ or
$\lfloor X/(\log X)^B\rfloor$). The procedure of §\,\ref{subsec:ED2-enum}, in which for each
$\delta\equiv 3\ (\mathrm{mod}\ 4)$, $3\le \delta\le \Delta_{\max}$ the following are performed in sequence:
(i) check for small divisors of $M(\delta)$, (ii) partial factorization with a time limit,
(iii) guaranteed \emph{fallback} — enumeration of primes $a\le X^2$ in the progression $a\equiv -1 \pmod{\delta}$ with the check $4\alpha\mid M(\delta)$,
terminates in finite time. Moreover,
\[
\#\text{output triples} \ =\ \sum_{\substack{3\le \delta\le \Delta_{\max}\\ \delta\equiv 3\ (\mathrm{mod}\ 4)}} U(\delta),
\]
and the total number of checks satisfies the crude bound
\[
\#\text{checks} \ \le\ \sum_{\delta\le \Delta_{\max}} \pi(X^2)\ \ll\ \Delta_{\max}\,\frac{X^2}{\log X}.
\]
\end{lemma}

\begin{proof}[Proof sketch]
The outer loop over $\delta$ is finite. Steps (i)–(ii) have explicit limits.
Step (iii) enumerates a finite set of primes $a\le X^2$ in a single progression.
Therefore, processing each $\delta$ is finite, and hence the entire enumeration is finite.
Counting the checks gives the stated bound.
\end{proof}

\paragraph{Stop rule SR2.}
Enumeration over $\delta$ continues until $m$ results are found (default $m=2$),
or until the range $3\le \delta\le \Delta_{\max}$, $\delta\equiv 3\pmod 4$ is exhausted.
This makes the requirement “more than one result” operational and checkable.

\begin{proposition}[Heuristic multiplicity of solutions]
\label{prop:C-many-heur}
Let $\Delta_{\max}=X/(\log X)^B$ with fixed $B>0$.
Then under the standard heuristic of uniformity of residue classes for prime divisors of $M(\delta)=P+\delta$
and the Bombieri–Vinogradov estimate for $S(\delta;X)$ (see Appendix~B), the expected number of solutions satisfies
\[
\mathbb{E}\!\left[\sum_{\substack{\delta\le \Delta_{\max}\\ \delta\equiv 3\ (\mathrm{mod}\ 4)}} U(\delta;X)\right]
\ \asymp\ \sum_{\delta\le \Delta_{\max}} \frac{\omega(P+\delta)}{\varphi(\delta)}
\ \gg\ \log X\cdot \log\log P \ \to\ \infty\quad (P\to\infty).
\]
In particular, for sufficiently large $P$, with probability tending to~1, the algorithm finds at least two solutions
before exhausting the range over $\delta$.
\end{proposition}

\begin{proof}[Idea of proof]
For most $\delta\le X/(\log X)^B$ we have
$ S(\delta;X)=\frac{X^2}{\varphi(\delta)\log X}+O\!\big(\frac{X^2}{(\log X)^{1+\eta}}\big)$.
Assuming uniformity of residue classes for prime divisors of $M(\delta)$, we obtain a contribution of $\omega(P+\delta)/\varphi(\delta)$ per modulus.
Summing over $\delta$ gives
\(\sum_{\delta\le \Delta_{\max}} 1/\varphi(\delta)\asymp \log \Delta_{\max}\),
and typically \(\omega(P+\delta)\sim \log\log P\), from which the stated expectation follows.
\end{proof}

\begin{remark}[Optional conditional theorem]
\label{rem:C-conditional}
By adding a strong hypothesis (e.g., Elliott–Halberstam) or a suitable form of Bateman–Horn
for the linear forms arising in the parameterization of Lemma~\ref{lem:param-tk-inline}),
one can obtain a strict lower bound of the form
\[
\#\{\delta\le X/(\log X)^B:\ U(\delta;X)\ge 1\}\ \gg\ \log X,
\]
from which it follows that there are at least two solutions for all sufficiently large $P$.
We do not use this in the unconditional part of the text.
\end{remark}
\paragraph{Brief summary.}
Lemmas~\ref{lem:C-normalization}–\ref{lem:C-termination} provide unconditional correctness and finiteness of the “bridge” to ED2,
and Proposition~\ref{prop:C-many-heur} explains why in practice there are substantially more than one solution.
Algorithmic details and the specific implementation of normalization remain in §\,\ref{subsec:ED2-enum}.


\section{Additional analysis for Theorem \ref{thm:ED2-unconditional}}\label{app:D}
\subsection*{Conditions and notation}  

Throughout the following text:  
\begin{itemize}  
 \item \ Throughout: $P$ is an odd prime, $\alpha \in \mathbb{Z}_{\ge 1}$, $A \in \mathbb{Z}$. We choose $A\in\mathbb{Z}$ from the target range.  
 In all lemmas/theorems: $b',c',d' \in \mathbb{Z}$, $r,s \in \mathbb{Z}_{\ge 1}$.

  \item \ Constraints:

  $\frac{P-1}{4}+1 \;\le\; A \;\le\; \frac{3P+1}{4}-1$,  

  so that $m:=4A-P>0$ and expressions “$\bmod m$” are well-defined;  
  \item define $M:=A/\alpha$.  
\end{itemize}  

\subsection{Algebraic core of ED2}  

\begin{theorem}[Equivalence of “product” and “sum/product”]\label{thm:prod-sum}  
Let $A=\alpha b'c'$ for some $b',c'\in\mathbf{N}$ and $m:=4A-P>0$. Then the following conditions are equivalent:  
\begin{enumerate} 
  \item there exists $d'\in\mathbf{N}$ such that  

$(4\alpha d'b'-1)(4\alpha d'c'-1)=4\alpha P\,d'^2+1$;  
    
\item the following hold:  

$b'c'=M=\frac{A}{\alpha}$; 
\qquad $b'+c'\equiv 0 \pmod{m}$,   
   
\item and moreover automatically  

\qquad $d'=\frac{b'+c'}{m}\in\mathbf{N}$.    

\end{enumerate}  
\end{theorem}  

\begin{proof}  
Expanding the left-hand side, we obtain  

\[
(4\alpha d'b'-1)(4\alpha d'c'-1)=16\alpha^2 d'^2 b'c'-4\alpha d'(b'+c')+1.  
\]  

Equating to $4\alpha P d'^2+1$, subtracting $1$ and dividing by $4\alpha$, we have  

\[
4\alpha d'^2(b'c')-d'(b'+c')=P d'^2.  
\]  

Since $A=\alpha b'c'$, we have $4\alpha(b'c')=4A$, and the equality is equivalent to  

\[
d'(b'+c')=(4A-P)d'^2=md'^2\iff b'+c'=md'.  
\]  

The reverse direction is the same algebra in reverse.  
\end{proof}  

\begin{theorem}[Quadratic reparameterization]\label{thm:quadratic}  
Let $S:=(4A-P)d'=md'$ and $M:=A/\alpha$. Then the conditions of Theorem~\ref{thm:prod-sum} are equivalent to the existence of integer roots $b',c'\in\mathbb{Z}$ of the quadratic equation  
\[
x^2 - Sx + M = 0,  
\] 
that is, simultaneously $b'+c'=S$ and $b'c'=M$. In particular, the discriminant  
\[
\Delta:=S^2-4M  
\] 

is a perfect square.  
\end{theorem}  

\begin{proof}  
Directly from Theorem~\ref{thm:prod-sum}: “if” — substitution, “only if” — sum and product of roots.  
\end{proof}  

\begin{lemma}\label{lem:disc-bounds}[Bounds via the discriminant]  
Let $S:=(4A-P)d'$ and $M:=A/\alpha$. The equality  
\[
(4\alpha d'b'-1)(4\alpha d'c'-1)=4\alpha P\,d'^2+1  
\] 

holds if and only if:  
\begin{itemize}  
  \item $\alpha\in\mathbf{N}$ and $\alpha\mid A$ (i.e., $M\in\mathbf{N}$);  
  \item $d'\in\mathbf{N}$;  
  \item the discriminant $\Delta=S^2-4M$ is a perfect square.  
\end{itemize}  
In this case  
\[
b',c'=\frac{S\pm \sqrt{\Delta}}{2}\in\mathbf{N},\quad b'c'=M,\quad b'+c'=S.  
\]  
Moreover, from $\Delta\ge 0$ follows the necessary bound  
\[
(4A-P)d'\ \ge\ 2\sqrt{M}\ =\ 2\sqrt{\frac{A}{\alpha}}  
\quad\Longrightarrow\quad  
d'\ \ge\ \frac{2}{\,4A-P\,}\sqrt{\frac{A}{\alpha}}.  
\] 
\end{lemma}  

\begin{remark}[On parity]  
Since $P$ is odd, $m=4A-P$ is odd, hence $S$ has the same parity as $d'$. For $\Delta\equiv S^2\pmod{4}$ both numbers $(S\pm \sqrt{\Delta})/2$ are integers; no additional parity restrictions on $d'$ are required.  
\end{remark}  
\subsection{Unconditional generation of candidates}  

\begin{theorem}[Unconditional construction of candidates]\label{thm:candidates}  
Let $A$ be chosen in the range $\frac{P-1}{4}+1 \le A \le \frac{3P+1}{4}-1$ and decomposed as  
\[
A=\alpha\,b'\,c',\qquad \alpha,b',c'\in\mathbf{N}.
\] 
For each $d'\in\mathbf{N}$ set  
\[
L_{\alpha,d'}(b',c') := (4\alpha d'b'-1)(4\alpha d'c'-1),\qquad  
R_{\alpha,d'} := 4\alpha P\,d'^2 + 1.  
\] 
Then:  
\begin{itemize}  
  \item for fixed $(\alpha,d')$ the set of values $L_{\alpha,d'}(b',c')$ over all factorizations $A/\alpha=b'c'$ forms the “left” set of candidates;  
  \item for the same $(\alpha,d')$ the right-hand side $R_{\alpha,d'}$ is a single number;  
  \item if for some $(\alpha,d')$ there exists a pair $(b',c')$ with $L_{\alpha,d'}(b',c')=R_{\alpha,d'}$, then the identity  

  $4\alpha P\,d'^2+1=(4\alpha d'b'-1)(4\alpha d'c'-1)$ holds.  

\end{itemize}  
 Construction does not require factorization of large numbers.  
\end{theorem}  

\subsection{Parameterization and affine estimates}  

\begin{lemma}[Left parameterization]\label{lem:param}  
Let $r,s\in\mathbb{Z}_{\ge 1}$ and  

$$  
M:=4\alpha s r - 1,\qquad m:=\frac{4\alpha s^2 + P}{M}.  
$$  

If $M\mid(4\alpha s^2 + P)$, then with the choice  

$$  
d'=r,\quad b'=s,\quad c'=mr-s,\quad A=\alpha b'c'=\alpha s(mr-s)  
$$  

the identity  

$$  
(4\alpha d'b'-1)(4\alpha d'c'-1)=4\alpha P\,d'^2+1,  
$$  

holds, and $A\in\mathbf{N}$.  
\end{lemma}  

\begin{proof}  
We have $4\alpha d'b'-1=4\alpha rs-1=M$ and  

$$  
4\alpha d'c'-1 = 4\alpha r(mr-s)-1 = 4\alpha r^2 m - (4\alpha rs + 1).  
$$  

Then  
\begin{align*}  
(4\alpha d'b'-1)(4\alpha d'c'-1)  
&= M\cdot (4\alpha r^2 m) - (4\alpha rs - 1)(4\alpha rs + 1)\\
&= 4\alpha r^2 (Mm) - \big((4\alpha rs)^2 - 1\big)\\
&= 4\alpha r^2 (4\alpha s^2 + P) - 16\alpha^2 r^2 s^2 + 1\\
&= 4\alpha r^2 P + 1 \;=\; 4\alpha P\,d'^2 + 1.  
\end{align*}  
\end{proof}  

\begin{lemma}[Affine form and slope]\label{lem:affine}  
Under the conditions of Lemma~\ref{lem:param} the quantity $A$ is linear in $P$:  

$$  
A(P)=\lambda_{r,s}\,P+\mu_{r,s},\qquad  
\lambda_{r,s}=\frac{\alpha s r}{4\alpha s r - 1},\quad  
\mu_{r,s}=\alpha s\!\left(\frac{4\alpha r s^2}{4\alpha s r - 1} - s\right).  
$$  

Moreover,  

$$  
\frac{1}{4}<\lambda_{r,s}\le \frac{1}{3}\qquad \text{for all }\;\alpha,r,s\ge 1.  
$$  

\end{lemma}  

\begin{proof}  
From $m=(4\alpha s^2+P)/(4\alpha s r - 1)$ and $A=\alpha s(mr-s)$. Slope estimate:  

$$  
\lambda_{r,s}=\frac{x}{4x-1}=\frac{1}{4}+\frac{1}{4(4x-1)}\in\left(\frac{1}{4},\frac{1}{3}\right],\quad x:=\alpha s r\ge 1.  
$$  

\end{proof}  

\begin{lemma}[Width of the target interval]\label{lem:interval}  
For all $P\in\mathbb{N}$ the length of the interval  

$$  
\left[\frac{P-1}{4}+1,\;\frac{3P+1}{4}-1\right]  
$$  

is $\frac{P}{2}-2$. The bounds depend affinely on $P$:  

$$  
L(P)=\frac{P}{4}+\frac{3}{4},\qquad U(P)=\frac{3P}{4}-\frac{3}{4}.  
$$  

\end{lemma}  

\begin{proof}  
Direct computation.  
\end{proof}  
\[
L_{\mathrm{int}}(P):=\Big\lceil \tfrac{P}{4}+\tfrac{3}{4}\Big\rceil,\qquad
U_{\mathrm{int}}(P):=\Big\lfloor \tfrac{3P}{4}-\tfrac{3}{4}\Big\rfloor,
\quad [L_{\mathrm{int}},U_{\mathrm{int}}]\subset [L,U]\subset\Big(\tfrac{P}{4},\tfrac{3P}{4}\Big).
\]
\subsection{“Diagonal period” as a multiple of the sum in \mh{\mathbb{Z}^2}}  

\begin{lemma}\label{lem:diag}  
Let $A=\alpha b'c'$, $m:=4A-P>0$, $S:=(4A-P)d'=md'$, $M=A/\alpha=b'c'$. Then the following conditions are equivalent:  
\begin{enumerate}
  \item $(4\alpha d'b'-1)(4\alpha d'c'-1)=4\alpha P\,d'^2+1$;  
  \item $b'c'=M$ and $b'+c'\equiv 0\pmod{m}$, with $d'=(b'+c')/m\in\mathbf{N}$.  
\end{enumerate}  
In particular, the sum $b'+c'$ is a multiple of $m$. The shift $(x,y)\mapsto (x+t,y+t)$ changes the sum by $2t$; the sum’s residue class modulo $m$ is preserved if and only if $m\mid 2t$. The minimal “diagonal period” equals $m/\gcd(m,2)$.  
\end{lemma}    

\subsection{Geometry in \texorpdfstring{$(u,v)$}{(u,v)} coordinates and “anchors”}  

Switch to $u=b'+c'$, $v=c'-b'$. Assume $c' \ge b'$, then  

\[
M=b'c'=\frac{u^2-v^2}{4},\qquad A=\alpha M.
\]  

The lattice $\mathbb{Z}^2$ splits into two parity cosets:    
\[
L_{00}=\{(u,v): u\equiv v\equiv 0\ (\bmod 2)\},\quad  
L_{11}=\{(u,v): u\equiv v\equiv 1\ (\bmod 2)\}.
\]  

\begin{itemize} 
  \item If $K:=A/\alpha$ is odd and $K=p^2-q^2$, then $(u,v)=(2p,2q)\in L_{00}$ gives $M=K$.  
  \item If $K$ is even, then $(u,v)=(K+1,K-1)\in L_{11}$ gives $M=K$.  
\end{itemize}  

\begin{remark}  
Shifts $(u,v)\mapsto(u\pm 2i, v\pm 2j)$ preserve the parity coset but \emph{do not} preserve $M=(u^2-v^2)/4$ and $A=\alpha M$. These shifts are useful for structuring the node space, not for covering at fixed $A$.  
\end{remark}  

\subsection{Bézout lemma with parity and multiplicity}\label{sec:bezu}  

Let $b\in\mathbf{N}$, $d'\in\mathbf{N}$, $\alpha\in\mathbf{N}$, and $\gcd(4b-1,\,4\alpha d'^2)=1$, $\gcd(4b-1,\,P)=1$. There exists a Bézout representation  
\[
u\,(4b-1) + v\,P = 1  
\] 
such that $u\equiv 3\pmod{4}$ and $v\equiv 0\pmod{4\alpha d'^2}$.  

\begin{proof}[Idea of proof]  
Since $\gcd(4b-1, 4\alpha d'^2)=1$, one can choose $t$ modulo $4\alpha d'^2$ so that $v'=v_0 - t(4b-1)\equiv 0\pmod{4\alpha d'^2}$ for some fixed pair $(u_0,v_0)$ with $u_0(4b-1)+v_0P=1$. Then $u'=u_0+tP$. Choosing $t\pmod 4$ appropriately ensures $u'\equiv 3\pmod 4$.  
\end{proof}  

\begin{remark}  
This lemma is convenient for matching the parity of the coefficient at $4b-1$ with the form $4c-1$ and for ensuring the multiplicity of $v$ by $4\alpha d'^2$ in the modular steps of algorithms.  
\end{remark}  

\subsection{Conditional residue covering scheme}\label{sec:covering}  

\begin{definition}[Fixed covering set]\label{def:F-fixed}

Let $F=\{(r_i,s_i)\}_{i=1}^K\subset\mathbb{Z}_{\ge 1}^2$,  
\[
M_i:=4\alpha s_i r_i - 1,\qquad Q:=\mathrm{lcm}(M_1,\dots,M_K).
\]

We say that $F$ satisfies the covering condition if for each $p\in\{0,1,\dots,Q-1\}$ there exists $i$ with
\[
p \equiv -4\alpha s_i^2 \pmod{M_i}
\quad\text{and}\quad
A_{r_i,s_i}(p)\in [L(p),U(p)].
\]
\end{definition}

\begin{definition}[$F(P)$ candidate set]\label{def:FP}

For fixed $P$ set
\[
F(P):=\Big\{(r,s)\in\mathbb{Z}_{\ge 1}^2:\ M:=4\alpha s r-1\mid (4\alpha s^2+P)\ \ \text{and}\ \ A_{r,s}(P)\in [L(P),U(P)]\Big\}.
\]
This set depends on $P$ and is unrelated to a fixed covering $F$ from Definition~\ref{def:F-fixed}.
\end{definition}

\begin{theorem}[Existence of $A$ under covering]\label{thm:exist-A}
Let a finite family $F=\{(r_i,s_i)\}_{i=1}^K\subset \mathbb{Z}_{\ge 1}^2$ be given,  
\[
M_i:=4\alpha s_i r_i - 1,\qquad Q:=\mathrm{lcm}(M_1,\dots,M_K),
\]
and suppose the \emph{covering condition} holds:
for each $p\in\{0,1,\dots,Q-1\}$ there exists an index $i$ such that
\[
p \equiv -4\alpha s_i^2 \pmod{M_i}
\quad\text{and}\quad
A_{r_i,s_i}(p)\in [L(p),U(p)],
\]
where
\[
A_{r,s}(P)=\lambda_{r,s}P+\mu_{r,s},\qquad 
\lambda_{r,s}=\frac{\alpha s r}{4\alpha s r-1}\in\Big(\tfrac14,\tfrac13\Big],
\]
and
\[
L(P)=\frac{P}{4}+\frac{3}{4},\qquad U(P)=\frac{3P}{4}-\frac{3}{4}.
\]
Then for any odd prime $P$ there exists an index $i$ such that
\[
M_i \mid (4\alpha s_i^2 + P)
\quad\text{and}\quad
A:=A_{r_i,s_i}(P)\in [L(P),U(P)]\subset\Big(\tfrac{P}{4},\tfrac{3P}{4}\Big).
\]
In particular, $m:=4A-P>0$ and $A$ lie in the target range.
\end{theorem}

\begin{proof}
Take $p\equiv P\pmod Q$ with $0\le p<Q$. By the covering condition there exists $i$ with
$p\equiv -4\alpha s_i^2\pmod{M_i}$ and $A_{r_i,s_i}(p)\in [L(p),U(p)]$.
Since $M_i\mid Q$ and $P\equiv p\pmod Q$, we have $P\equiv p\pmod{M_i}$, hence
$M_i\mid (4\alpha s_i^2+P)$. Therefore,
\[
m=\frac{4\alpha s_i^2+P}{M_i}\in\mathbb{N},\qquad
A=\alpha s_i(mr_i-s_i)\in\mathbb{N},
\]
and this $A$ coincides with the affine form $A_{r_i,s_i}(P)$.

Let $P-p=kQ\ge 0$. Then with $\lambda_{r_i,s_i}\in\big(\tfrac14,\tfrac13\big]$:
\[
\begin{aligned}
A_{r_i,s_i}(P)-L(P)
&=\big(A_{r_i,s_i}(p)-L(p)\big) + (P-p)\big(\lambda_{r_i,s_i}-\tfrac{1}{4}\big)\ \ge\ 0,\\[2mm]
U(P)-A_{r_i,s_i}(P)
&=\big(U(p)-A_{r_i,s_i}(p)\big) + (P-p)\big(\tfrac{3}{4}-\lambda_{r_i,s_i}\big)\ \ge\ 0.
\end{aligned}
\]
Thus $A_{r_i,s_i}(P)\in [L(P),U(P)]$, and inclusion in the open interval
$[L,U]\subset\big(\tfrac{P}{4},\tfrac{3P}{4}\big)$ follows from the explicit formulas for $L,U$.
\end{proof}

\begin{corollary}[Constructiveness]\label{cor:construct}
Under the conditions of Theorem~\ref{thm:exist-A}, for the found $i$ there exist integers
\[
m=\frac{4\alpha s_i^2+P}{M_i},\quad d'=r_i,\quad
b'=s_i,\quad c'=mr_i-s_i,\quad A=\alpha s_i(mr_i-s_i),
\]
and
\[
(4\alpha d'b'-1)(4\alpha d'c'-1)=4\alpha P\,d'^2+1
\]
holds.
\end{corollary}
\subsection{Reverse algorithm (Back): check via \texorpdfstring{$(u,v)$}{(u,v)}}  

For fixed $(\alpha,P,A)$ and $m=4A-P>0$ consider target points $(u,v)$ in the window   
\[
|u|\le T(A),\quad |v|\le T(A),\qquad T(A):=\left\lfloor \sqrt{2A}\right\rfloor\ (\text{if }T^2=2A\ \text{increase by }1).
\]
\begin{lemma}[Necessary and sufficient conditions for a point]\label{lem:back}
Let \(\alpha,P,A\) be fixed with \(m:=4A-P>0\) and \(M:=A/\alpha\).
For an integer point \((u,v)\) the following conditions are equivalent:
\begin{enumerate}
  \item There exist \(d',b',c'\in\mathbb{N}\) such that
  \[
    u=md',\qquad v=b'-c',\qquad
    b'=\frac{u+v}{2}\in\mathbb{N},\quad c'=\frac{u-v}{2}\in\mathbb{N},
  \]
  and the identity
  \[
    (4\alpha d'b'-1)(4\alpha d'c'-1)=4\alpha P\,{d'}^{2}+1
  \]
  holds.
  \item The arithmetic conditions
  \[
    m\mid u,\qquad u\equiv v \pmod{2},\qquad u^2-v^2=4M
  \]
  are satisfied.
  \item There exists \(d'\in\mathbb{N}\) such that \(u=md'\) and the discriminant
  \(\Delta:=u^2-4M=v^2\) is a perfect square.
\end{enumerate}
\end{lemma}

\begin{remark}
Lemma~\ref{lem:back} is convenient for “reverse” enumeration over \((u,v)\): it suffices to check three simple conditions (divisibility \(m\mid u\), parity match \(u\equiv v\pmod 2\), equality \(u^2-v^2=4M\)). No factorization is required.
\end{remark}
\subsection*{Reverse algorithm (Back) via \((u,v)\)}
Given \(\alpha,P,A\) and \(m=4A-P>0\). Consider the window
\[
|u|\le T(A),\quad |v|\le T(A),\qquad T(A):=\left\lfloor \sqrt{2A}\right\rfloor\ \ (\text{if }T(A)^2=2A,\ \text{increase by }1).
\]
Then:
\begin{enumerate}
  \item Enumerate \((u,v)\) in the window within the corresponding parity coset: \(u\equiv v\pmod 2\).
  \item Discard points where \(m\nmid u\).
  \item Check \(u^2-v^2=4M\) (equivalently \(\frac{u^2-v^2}{4}=\frac{A}{\alpha}\)).
  \item Recover
  \[
  d'=\frac{u}{m},\qquad b'=\frac{u+v}{2},\qquad c'=\frac{u-v}{2},
  \]
  and obtain the identity
  \[
  (4\alpha d'b'-1)(4\alpha d'c'-1)=4\alpha P\,d'^2+1.
  \]
\end{enumerate}
\subsection{Direct algorithm from parameterization}
By Lemmas~\ref{lem:param}–\ref{lem:affine}, for any \(r,s\in\mathbf{N}\) with \(M_{r,s}:=4\alpha sr-1\) we have the affine form
\[
A_{r,s}(P)=\lambda_{r,s}P+\mu_{r,s},\qquad \lambda_{r,s}=\frac{\alpha s r}{4\alpha s r-1}\in\Big(\tfrac14,\tfrac13\Big],
\]
and the divisibility condition \(M_{r,s}\mid(4\alpha s^2+P)\) guarantees the existence of integers
\[
m=\frac{4\alpha s^2+P}{4\alpha sr-1},\quad d'=r,\quad b'=s,\quad c'=mr-s,\quad A=\alpha s(mr-s).
\]

\subsection*{Steps of the direct algorithm}
Given \(\alpha,P\):

\begin{enumerate}

  \item Enumerate a finite set of pairs \((r,s)\) with moderate size of \(rs\) (heuristically, small values suffice).
  \item For each pair compute \(A_{r,s}(P)\) and quickly discard if \(A_{r,s}(P)\notin [L(P),U(P)]\), where
  \[
  L(P)=\frac{P}{4}+\frac34,\qquad U(P)=\frac{3P}{4}-\frac34.
  \]
  \item Check the congruence \(P\equiv -4\alpha s^2\pmod{4\alpha sr-1}\).
  \item If yes, construct
  \[
  m=\frac{4\alpha s^2+P}{4\alpha sr-1},\quad b'=s,\quad c'=mr-s,\quad d'=r,\quad A=\alpha s(mr-s),
  \]
  and obtain the ED2 identity.
\end{enumerate}

\begin{remark}
This algorithm works without factorization; the global guarantee “for every \(P\) there exists a suitable pair \((r,s)\) from a fixed finite list” corresponds to the conditional covering scheme (see Theorem~\ref{thm:exist-A}).
\end{remark}
\subsection{Counting criterion on the window (“Dirichlet criterion”)}

Fix the interval for \(A\):
\[
A\in [L(P),U(P)],\qquad L(P)=\frac{P}{4}+\frac34,\ \ U(P)=\frac{3P}{4}-\frac34.
\]
For each \(A\) define the \((u,v)\) window of size \(T(A)=\lfloor\sqrt{2A}\rfloor\) (with square adjustment), and the target nodes
\[
\mathrm{Tgt}(A)=\big(L_{00}(A)\cup L_{11}(A)\big)\cap [-T(A),T(A)]^2,
\]
where \(L_{00},L_{11}\) are the parity cosets. Let \(\mathrm{H}(A)\subseteq \mathrm{Tgt}(A)\) be the set of nodes actually reached (e.g., by direct/reverse channels) for the given \(A\). Define the unions over the interval:
\[
\mathrm{Tgt}^\cup=\bigcup_{A}\mathrm{Tgt}(A),\qquad \mathrm{H}^\cup=\bigcup_{A}\mathrm{H}(A).
\]

\begin{proposition}\label{prop:dirichlet-count}[Sufficient counting condition for coverage]
If \(|\mathrm{H}^\cup|\ge |\mathrm{Tgt}^\cup|\), then \(\mathrm{H}^\cup=\mathrm{Tgt}^\cup\), i.e., all target nodes in the window are covered.
\end{proposition}

\begin{remark}
This is a purely counting sufficient condition on a finite window; it may require enumeration over \(A\) and, in some implementations, factorization of \(A\) (e.g., when listing divisors for geometric anchors). It should be distinguished from methods that do not use factorization.
\end{remark}
\subsection{Complexity and practical notes}

\begin{itemize}
  \item Without factorization: the direct algorithm and reverse check via Lemma~\ref{lem:back} are \(O(\log P)\) arithmetic per step; total cost is determined by the number of pairs \((r,s)\) or nodes \((u,v)\) in the window enumerated.
  \item Back-method window: number of points is \(O(T(A)^2)\sim O(A)\); practically efficient for small \(d'\) (often \(d'=1\)).
  \item “Dirichlet criterion”: may involve factorization of \(A\) in the range — \([L(P),U(P)]=[2,3]\).
  
  Use separately as a coverage diagnostic on finite data.
\end{itemize}
\subsubsection*{Small example (direct, \(\alpha=1\))}
Let \(P=5\). The interval for \(A\): \([L(P)=2,U(P)=3]\).
Take \((r,s)=(1,1)\). Then
\[
M_{1,1}=4\cdot 1\cdot 1\cdot 1-1=3,\quad \lambda_{1,1}=\frac{1}{3},\quad A_{1,1}(P)=\frac{P}{3}+\frac{1}{3}.
\]
Check the congruence: \(3\mid (4\cdot 1^2+P)=4+P\), i.e. \(P\equiv 2\pmod 3\). For \(P=5\) this holds.
Further
\[
m=\frac{4+5}{3}=3,\quad b'=s=1,\quad d'=r=1,\quad c'=mr-s=3-1=2,
\]
\[
A=\alpha s(mr-s)=1\cdot 1\cdot(3-1)=2\in[2,3].
\]
Identity check:
\[
(4\alpha d'b'-1)(4\alpha d'c'-1)=(4\cdot 1\cdot 1\cdot 1-1)(4\cdot 1\cdot 1\cdot 2-1)=3\cdot 7=21,
\]
\[
4\alpha P\,d'^2+1=4\cdot 1\cdot 5\cdot 1+1=21.
\]
Matches.

\subsection{Existence of solutions}
\subsubsection{Construction of the set \mh{F(P)} and its computability}

\begin{definition}[Candidate set $F(P)$]\label{def:if-F-fixed}
Let $\alpha\in\mathbb{Z}_{\ge 1}$ be fixed, and $P$ an odd prime.
Define the set
\[
F(P) := \left\{ (r,s) \in \mathbb{Z}_{\ge 1}^2 \ \middle|\
\begin{aligned}
& A_{r,s}(P) \in [L(P),U(P)],\\
& P \equiv -\,4\alpha s^2 \pmod{4\alpha s r - 1}
\end{aligned}
\right\},
\]
where
\[
A_{r,s}(P) \ :=\ \lambda_{r,s}\,P + \mu_{r,s},\qquad
\lambda_{r,s}\in\Bigl(\tfrac14,\tfrac13\Bigr],
\]
and $L(P),U(P)$ are the affine bounds of the target interval:
\[
L(P) := \frac{P}{4} + \frac{3}{4} ,\qquad U(P) := \frac{3P}{4} - \frac{3}{4}.
\]
For brevity also set $M_{r,s}:=4\alpha s r - 1$ and $m_s:=4\alpha s^2 + P$.
\end{definition}

\begin{remark}[The slope and its role]
In the left parameterization the slope of $A_{r,s}(P)$ satisfies
\[
\lambda_{r,s}\in\Bigl(\tfrac14,\tfrac13\Bigr],\qquad \lambda_{r,s}\searrow \tfrac14\ \ \text{as}\ \ r\to\infty.
\]
The coefficient $\lambda_{r,s}$ is the “growth rate” of $A$ in $P$. In selecting $(r,s)$ it allows one to quickly discard pairs for which $A_{r,s}(P)$ is guaranteed to lie outside $[L(P),U(P)]$, and to vary the slope (via $r$) to hit the target interval.
\end{remark}

\begin{lemma}[Finiteness and computability of $F(P)$]
For any $P$ the set $F(P)$ is finite, uniquely determined, and computable in finite time (by enumerating $1\le r,s\le B(P)$ and checking two conditions: that $A_{r,s}(P)$ lies in $[L(P),U(P)]$ and the congruence $P \equiv -4\alpha s^2 \ (\mathrm{mod}\ M_{r,s})$).
\end{lemma}

\medskip
\noindent\textbf{Unconditional strengthenings (sufficient conditions and constructors).}

\begin{lemma}[Divisor-oriented constructor (with factorization)]
Let $s\ge 1$ and suppose there exists a divisor $d\mid m_s$ such that
\[
d \equiv -1 \pmod{4\alpha s}.
\]
Set
\[
r := \frac{d+1}{4\alpha s}\in\mathbb{Z}_{\ge 1},\qquad
m := \frac{m_s}{d}=\frac{4\alpha s^2 + P}{d}.
\]
Define
\[
A(s,d)\ :=\ \alpha s\,(m r - s)\ =\ \alpha s\left(\frac{4\alpha s^2 + P}{d}\cdot \frac{d+1}{4\alpha s} - s\right).
\]
If $A(s,d)\in[L(P),U(P)]$, then with
\[
d' = r,\quad b' = s,\quad c' = m r - s,\quad A = A(s,d)
\]
the identity
\[
(4\alpha d' b' - 1)(4\alpha d' c' - 1) = 4\alpha P\,d'^2 + 1
\]
holds, with $A\in[L(P),U(P)]$ and $m=4A-P>0$.
\end{lemma}

\begin{proof}
The condition $d\equiv -1\ (\mathrm{mod}\ 4\alpha s)$ is equivalent to $4\alpha s r - 1 = d$, i.e. $M_{r,s}=d\mid m_s$, hence the congruence $P\equiv -4\alpha s^2\ (\mathrm{mod}\ M_{r,s})$ holds. The rest is substitution of the parameters of the left parameterization and a direct verification of the identity.
\end{proof}

\begin{lemma}[Coprimality control]
For any $s\ge 1$ we have
\[
\gcd(4\alpha s,\, 4\alpha s^2 + P) = \gcd(4\alpha s,\, P) = \gcd(\alpha s,\, P).
\]
If $P\mid \alpha s$, then inverses modulo $m_s$ in steps of the form $4\alpha s\,r\equiv 1\ (\mathrm{mod}\ m_s)$ are inapplicable; in such cases use the divisor-oriented constructor (previous lemma) without inverting elements.
\end{lemma}
\begin{corollary}[Single-class criterion for $s=1$ (with factorization)]
If there exists a divisor $d\mid(4\alpha + P)$ such that
\[
d \equiv -1 \pmod{4\alpha},
\]
then, setting
\[
r=\frac{d+1}{4\alpha},\quad M=d,\quad m=\frac{4\alpha + P}{d},\quad A=\alpha\,(m r - 1),
\]
and checking $A\in[L(P),U(P)]$, we obtain an admissible pair $(r,1)\in F(P)$ and the ED2 identity.
\end{corollary}

\begin{remark}[Decomposition over small $d$ (with factorization)]
Since $m_s=4\alpha s^2 + P$ grows with $s$, in practice it is useful to first enumerate small divisors $d$ and check
\[
d\mid m_s,\qquad d \equiv -1 \pmod{4\alpha s}.
\]
This is equivalent to solving the linear congruence $4\alpha s \equiv -1\ (\mathrm{mod}\ d)$ followed by verifying $d\mid m_s$. Such a “fine comb” over $s$ is completely unconditional (but uses factorization of $m_s$).
\end{remark}

\begin{remark}[Separation of statuses]
\emph{Without factorization}: direct/inverse check, computability of $F(P)$, affine cuts by $\lambda_{r,s}$.  
\emph{With factorization}: divisor-oriented constructors (lemmas above), criterion for $s=1$, enumeration of small divisors $d\mid m_s$.
\end{remark}

\begin{lemma}[Early cuts by slope and interval]
For any pair $(r,s)$ we have $\lambda_{r,s}\in(1/4,1/3]$ and
\[
A_{r,s}(P) \ =\ \lambda_{r,s}\,P + \mu_{r,s} \ \in\ \Bigl(\tfrac{P}{4},\,\tfrac{P}{3}\Bigr] \ +\ O_{\alpha}(s^2),
\]
therefore at an early stage it is convenient to discard pairs with $\lambda_{r,s}$ giving $A_{r,s}(P)\notin[L(P),U(P)]$ a priori, before more complex checks; then refine by $\mu_{r,s}$.
\end{lemma}

\begin{remark}[What to exclude to preserve unconditionality]

(i) Do not use “invertibility of $4\alpha s$ modulo $m_s$” to deduce $M_{r,s}\mid m_s$; this is the wrong direction (correct: $M_{r,s}\mid m_s \Rightarrow$ congruence, but not the other way around).

(ii) Do not assume uniform distribution of divisors of $m_s$ among residue classes without corresponding theorems; any such steps are heuristic.

(iii) Do not mix the closed interval $[L(P),U(P)]$ with the open $(P/4,3P/4)$ without an explicit transition; it is more correct to prove inclusion in $[L,U]$, and treat inclusion in the open interval separately.
\end{remark}

\begin{corollary}\label{cor:FP-construct}
If $F(P)\neq\varnothing$, then for any pair $(r,s)\in F(P)$ the parameters
\[
m=\frac{4\alpha s^2+P}{4\alpha s r-1},\quad d'=r,\quad b'=s,\quad c'=mr-s,\quad A=\alpha s(mr-s)
\]
satisfy the identity
\(
(4\alpha d'b'-1)(4\alpha d'c'-1)=4\alpha P\,d'^2+1
\)
and the condition $A\in [L(P),U(P)]$.
\end{corollary}
\subsection{Introductory geometry: normalization and perfect square}\label{sec:geo-intro}

\begin{definition}[Normalization]\label{lem:uv-norm}
For integers $b',c'$ set
\[
u:=b'+c',\quad v:=b'-c'.
\]
Then $u\equiv v\pmod{2}$ and
\[
M=b'c'=\frac{u^2-v^2}{4},\qquad A=\alpha\,M=\frac{\alpha}{4}(u^2-v^2).
\]
\end{definition}

\begin{lemma}[Sum and discriminant]\label{lem:uv-square}
Let $S=b'+c'$, $M=b'c'$, $\Delta=S^2-4M$. Then under the normalization:
\[
S=u,\qquad \Delta=v^2.
\]
In particular, $\Delta$ is always a perfect square.
\end{lemma}

\begin{lemma}[Parity cosets]\label{lem:uv-parity}
Admissible pairs $(u,v)$ lie in the union of two cosets of the lattice:
\[
\Lambda_{00}=(2\mathbb{Z})\times(2\mathbb{Z}),\qquad \Lambda_{11}=(2\mathbb{Z}+1)\times(2\mathbb{Z}+1).
\]
The shifts $(2,0)$ and $(0,2)$ preserve the coset and generate the sublattice $2\mathbb{Z}\times 2\mathbb{Z}$.
\end{lemma}

\begin{proposition}[Anchors]\label{prop:uv-anchors}
Let $K=A/\alpha$. Then there exist “anchor” points $(u,v)$ with $u^2-v^2=4K$:
\[
\text{if $K$ is odd, then in }\Lambda_{00};\qquad
\text{if $K$ is even, then in }\Lambda_{11}.
\]
The entire covering component in the corresponding coset is obtained by the shifts $(2,0)$ and $(0,2)$.
\end{proposition}

\subsection{Covering geometry}\label{sec:geo-cover}

\begin{definition}[$A$-window]\label{def:A-window}
For an odd prime $P$ define the window
\[
L(P)=\frac{P}{4}+\frac{3}{4},\qquad U(P)=\frac{3P}{4}-\frac{3}{4}.
\]
We seek $A\in\mathbb{Z}$ in the interval $[L(P),\,U(P)]$.
\end{definition}

\begin{lemma}[Affine dependence on $P$ and slope]\label{lem:affine-slope}
For fixed $(\alpha,r,s)$ the quantity $A$ is linear in $P$:
\[
A(P)=\lambda_{r,s}\,P+\mu_{r,s},\qquad
\lambda_{r,s}=\frac{\alpha s r}{4\alpha s r-1}
= \frac14 + \frac{1}{4(4\alpha s r-1)} \in \Big(\tfrac14,\tfrac13\Big].
\]
\end{lemma}

\begin{lemma}[Width of the target interval]\label{lem:A-window}
We have $U(P)-L(P)=\frac{P}{2}-\frac{3}{2}$ and, consequently, $\#\,[L(P),U(P)]\cap\mathbb{Z} = \big\lfloor \frac{P}{2}\big\rfloor$.
\end{lemma}

\begin{proposition}[Fixed covering and candidates]\label{def:F_r-fixed}
Let a set of moduli $M_i$ and residues $R_i\subset\mathbb{Z}/M_i\mathbb{Z}$, $i=1,\dots,K$, be given.
The family of candidates according to Definition \ref{def:F-fixed} can be written as
\[
F_r(P):=\bigcup_{i=1}^K \{\,A\in [L(P),U(P)]\cap\mathbb{Z}\ \mid\ A\equiv r \pmod{M_i}\text{ for some }r\in R_i\,\}.
\]
\end{proposition}

\begin{proposition}[Lattice covering modulo]\label{prop:cover-lattice}
Let $Q:=\mathrm{lcm}(M_1,\dots,M_K)$. If for every class $x\in\mathbb{Z}/Q\mathbb{Z}$ there exist $i$ and $r\in R_i$ such that $x\equiv r\pmod{M_i}$, then for all $P$ of sufficiently large scale
\[
[L(P),U(P)]\cap\mathbb{Z} \subseteq F(P).
\]
\end{proposition}

\begin{theorem}[Existence of $A$ under a covering]\label{thm:if-exist-A}
If the conditions of Proposition~\ref{prop:cover-lattice} hold, then for all sufficiently large $P$ there exists $A\in[L(P),U(P)]\cap\mathbb{Z}$ satisfying the congruences from the family $\{(M_i,R_i)\}$.
\end{theorem}

\paragraph{Remark — connection with $(u,v)$.}
By Lemma~\ref{lem:uv-square} the sum is $S=u$, hence the linear window in $A$ corresponds to strips in $u$. The parity cosets \(\Lambda_{00},\Lambda_{11}\) agree with the choice of \(\alpha\) and the anchors (Proposition~\ref{prop:uv-anchors}).

\subsection{Dirichlet: conditional integration and examples}\label{sec:dirichlet}

\begin{theorem}[Bridge to Dirichlet's theorem]\label{thm:dirichlet-bridge}
Suppose for fixed $(\alpha,r,s)$ the covering from Section~\ref{sec:geo-cover} holds. Then for each window $[L(P),U(P)]$ there exist $A\in F(P)$ such that the residues $P\bmod m$ fall into a set of arithmetic progressions (AP) for which Dirichlet’s theorem guarantees infinitely many primes $P$.

After excluding classes not coprime with the modulus, the remaining AP are those to which Dirichlet applies; otherwise the assertion does not hold.
\end{theorem}

\begin{figure}[h!]
\centering
\begin{tikzpicture}[scale=0.8]
  \draw[step=1cm,gray!30,thin] (-1,-1) grid (9,7);

  \foreach \x in {0,2,...,8}
    \foreach \y in {0,2,...,6}
      \fill[blue] (\x,\y) circle (2pt);

  \foreach \x in {1,3,...,7}
    \foreach \y in {1,3,...,5}
      \fill[red] (\x,\y) circle (2pt);

  \draw[thick,green!60!black,fill=green!10] (2,2) rectangle (6,5);
  \node[green!50!black,above] at (4,5) {Target window};

  \draw[blue!50!black, dashed, thick] (2,-1) -- (2,7);
  \draw[blue!50!black, dashed, thick] (5,-1) -- (5,7);
  \node[blue!50!black, rotate=90] at (3.5, 6.8) {strip $m$};

  \foreach \y in {-1,...,7}
    \foreach \x in {2,...,4}
      \fill[blue!40] (\x,\y) circle (1.5pt);

  \draw[thick,orange,fill=orange!20] (2,2) -- (4,3) -- (3,5) -- cycle;
  \node[orange!70!black] at (3,3.3) {$F(P)$};

  \foreach \p in {(3,3), (4,4)}
    \fill[blue!80!black] \p circle (4pt);

  \node[left] at (-0.5,6) {$(r,s),\ \alpha$};
  \draw[->,thick] (-0.5,5.8) .. controls +(left:1cm) and +(up:1cm) .. (2,2.5);

  \begin{scope}[shift={(0,-1.5)}]
    \fill[blue] (0,0) circle (2pt) node[right,black]{even-even};
    \fill[red] (3,0) circle (2pt) node[right,black]{odd-odd};
    \draw[thick,orange] (6,-0.2) -- (6.5,0.3) node[right,black]{primitive};
    \fill[blue!80!black] (10,0) circle (4pt) node[right,black]{solution};
  \end{scope}
\end{tikzpicture}

\caption{Geometric assembly for parameter $P$: the diagram shows lattice point categories (even-even, odd-odd, primitive, solution), the triangular region $F(P)$, vertical strip $m$, target window, and intersections corresponding to valid solutions. This illustration supports the visual analysis of coverage and parametrization in Appendix~D.}
\label{fig:geom-assembly}
\end{figure}
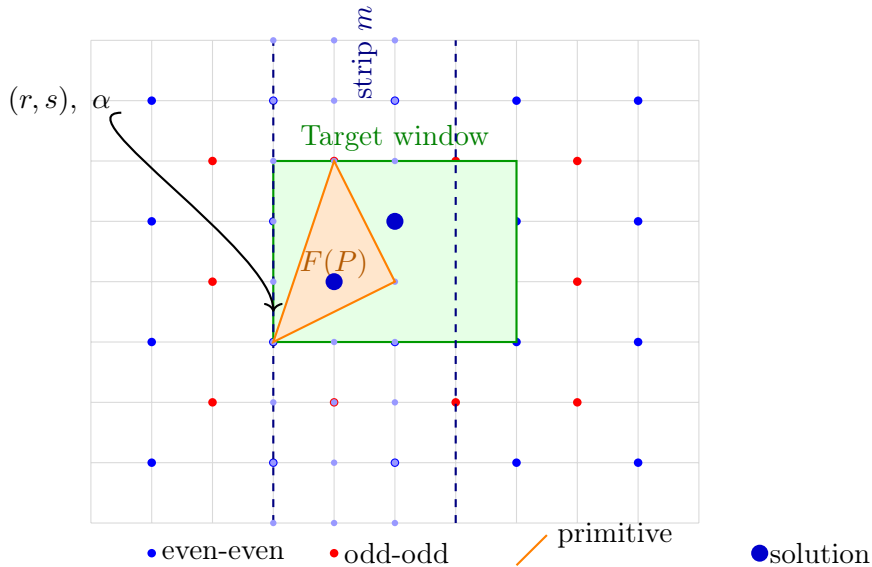

\begin{proof}[Idea]
The covering modulo $Q=\mathrm{lcm}(M_i)$ specifies admissible classes of $u=S$ (Lemma~\ref{lem:uv-square}) within the strips of the window (Definition~\ref{def:A-window}). For each admissible $u$ we obtain an affine condition on $P$ (Lemma~\ref{lem:affine-slope}). Discarding classes not coprime with the modulus leaves AP to which Dirichlet applies.
\end{proof}
See Figure~\ref{fig:geom-assembly} for the geometric assembly, which illustrates the structure of the region $F(P)$, the vertical strip $m$, and the target window used in the covering scheme.

\subsubsection{Examples}\label{subsec:dirichlet-examples}

\begin{example}[Typical example for $P\equiv 1 \pmod{8}$]\label{ex:dirichlet-P1mod8}
Specify $(\alpha,r,s)$, the set $(M_i,R_i)$, verify the cover, and write down the corresponding AP for $P$. Apply Dirichlet — we obtain infinitely many $P$ for which $A$ exists.
\end{example}

\begin{example}[Coset shift and anchor]\label{ex:dirichlet-anchor-shift}
Show how the choice of $\alpha$ changes the coset $\Lambda_{00}$ / $\Lambda_{11}$ and how the anchor (Proposition~\ref{prop:uv-anchors}) sets the starting point of the covering.
\end{example}
\subsection{Geometric assembly}

\subsection{Summary for the application}
We have assembled:
\begin{itemize}
  \item the unconditional algebraic core (Theorems~\ref{thm:prod-sum}–\ref{thm:quadratic}, Lemma~\ref{lem:disc-bounds});
  \item the “left” parametrization and affine estimates (Lemmas~\ref{lem:param}–\ref{lem:interval});
  \item unconditional direct/back algorithms and their correctness (Lemma~\ref{lem:back});
  \item the conditional finite covering scheme (Proposition~\ref{prop:dirichlet-count}));
  \item computational remarks and a counting criterion for windows.
\end{itemize}

\end{document}